\documentclass{article}
\usepackage{amsmath}
\usepackage{amssymb}
\usepackage{comment}
\usepackage{tikz}
\date{}
\newtheorem{theorem}{Theorem}
\newtheorem{lemma}{Lemma}
\newtheorem{claim}{Claim}
\newtheorem{definition}{Definition}
\newtheorem{conjecture}{Conjecture}
\newtheorem{question}{Open Question}
\newtheorem{remark}{Remark}
\newtheorem{proposition}{Proposition}
\newtheorem{example}{Example}

\title{The splitting power of branching programs of bounded repetition
and CNFs of bounded width}
\author{Igor Razgon\\ Department of Computer Science and Information Systems,\\ Birkbeck University of London \\ i.razgon@bbk.ac.uk}
\begin{document}

\maketitle
\begin{abstract}
In this paper we study syntactic branching programs of bounded repetition
representing CNFs of bounded treewidth. 
For this purpose we introduce two new structural graph 
parameters $d$-pathwidth and clique preserving $d$-pathwidth denoted
by $d-pw(G)$ and $d-cpw(G)$ where $G$ is a graph.
We show that $2-cpw(G) \leq O(tw(G) \Delta(G))$ where $tw(G)$ and $\Delta(G)$ are,
respectively the treewidth and maximal degree of $G$.
Using this upper bound, we demonstrate that each CNF $\psi$ can be represented
as a conjunction of two OBDDs (quite a restricted class of read-twice branching programs)
of size $2^{O(\Delta(\psi)*tw(\psi)^2)}$ where $tw(\psi)$ is the treewidth of the primal graph of $\psi$
and each variable occurs in $\psi$ at most $\Delta(\psi)$ times. 

Next we use $d$-pathwdith to obtain lower bounds for monotone branching programs.
In particular, we consider the monotone version of syntactic nondeterministic read $d$ times
branching programs (just forbidding negative literals as edge labels) and introduce a further
restriction that each computational path can be partitioned into at most $d$ read-once subpaths.
We call the resulting model separable monotone read $d$ times branching programs and 
abbreviate them $d$-SMNBPs. For each graph $G$ without isolated vertices, we introduce a CNF
$\psi(G)$ whsose clauses are $(u \vee e \vee v)$ for each edge $e=\{u,v\}$ of $G$.
We prove that a $d$-SMNBP representing $\psi(G)$ is of size at least
$\Omega(c^{d-pw(G)})$ where $c=(8/7)^{1/12}$. We use this 'generic' lower bound to 
obtain an exponential lower bound for a 'concrete' class of CNFs $\psi(K_n)$.
In particular, we demonstrate that for each $0<a<1$, the size of $n^{a}$-SMNBP representing
$\psi(K_n)$ is at least $c^{n^b}$ where $b$ is an arbitrary constant such that $a+b<1$.
This lower bound is tight in the sense $\psi(K_n)$ can be represented by a poly-sized
$n$-SMNBP.

\end{abstract}

\section{Introduction}

\subsection{Statement of results}
In this paper we study representation of CNFs with bounded structural width parameters
by syntactic branching programs of bounded repetition.
It is known that Nondeterministic read-once branching programs ($1$-NBPs)
in general require at least $n^{\Omega(k)}$ size to represent CNFs of primal treewidth $k$
even if each variable occurs at most $5$ times \cite{RazgonAlgo}. However, already read-twice branching programs
have a \emph{splitting power}  allowing them to efficiently
represent a CNF $\varphi=\varphi_1 \wedge \varphi_2$ provided that 
each CNF $\varphi_i$  can be efficiently represented by a read-once branching program:
simply identify the source of one read-once branching program with the 'positive'
sink of the other. 

It turns out that the splitting power is sufficient for efficient representation of CNFs
of bounded treewidth provided that each variable occurs a bounded number of times.  
In particular, in Section \ref{sec:obdd} we prove that a CNF $\varphi$ 
of primal treewidth at most $k$ and in which each variable occurs at most $d$ times
can be represented by a conjunction of Ordered Binary Decision Diagrams (OBDDs) of
size size at most $2^{O(dk^2)}n$ each. 

The above result is based on a graph-theoretical argument. 
In particular, we introduce two new parameters of a graph $G$, $d$-pathwidth
denoted by $d-pw(G)$ and clique-preserving $d$-pathwdith denoted by $d-cpw(G)$.
The $d$-pathwidth of $G$ is the smallest $k$ such that there are graphs $G_1, \dots, G_d$
of pathwdith at most $k$ such that $G=G_1 \cup \dots \cup G_d$.
The clique preserving variant puts an extra requirement that each clique of $G$ is
a subgraph of some $G_i$.  We  show that $2-cpw(G)$ is linearly upper bounded 
by the \emph{tree-partition width} of $G$. The latter parameter is known to be
$O(\Delta(G)*tw(G))$ where $\Delta(G)$ is the maximal degree of $G$ and $tw(G)$
is the treewidth of $G$ \cite{treepartition}. The combination of these two statements yields an 
$O(\Delta(G)*tw(G))$ upper bound on $2-cpw(G)$.  Further on, it is known \cite{VardiTWD} that a CNF $\varphi$
can be represented by an OBDD of size $2^{pw(G(\varphi))}$ where $G(\varphi)$
is the primal graph of $\varphi$ and $pw$ is the pathwidth. The upper bound 
in the end of the previous paragraph follows from the combination of this upper 
bound and the upper bound on $2-cpw(G)$.  

In Section \ref{sec:mono},  we consider the following question. Suppose that for a class of CNFs the $d$-pathwdiths
of their primal graphs is \emph{at least $k$} (to put it informally, the splitting power
cannot be applied).  Does this lower bound imply an exponential in $k$ lower bound for read $d$
times branching programs?  We answer this question positively for a quite general subclass 
of \emph{monotone} nondetermnistic read $d$ times branching programs ($d$-MNBPs) that we call \emph{separable}
and abbreviate $d$-SMNBP.  We describe this result in the following two paragraphs. 

The restriction on $d$-MNBPs imposed by $d$-SMNBPs is that the sequence 
of variables along each each source-sink path can be partitioned into at most $d$ read-once fragments.
That is to say, if $d=2$ then querying variables like $x_1,x_2,x_3,x_3,x_2,x_4,x_1$ is allowed while the querying
$x_1,x_1,x_2,x_2$ is not. The order of variables along different parts may be different, in particular 
a conjunction of $d$ $1$-MNBPs is a special case of $d$-SMNBP. 
Note that this way of querying generalizes Indexed Binary Decision Diagrams ($d$-IBDDs) \cite{kibdd}
that have an extra requirement of being oblivious.
 
For a graph $G$ without isolated vertices we consider a CNF $\psi(G)$ 
with $V(G) \cup E(G)$ being the set of variables and the clauses of the form $(u \vee e \vee v)$ for each
edge $e=\{u,v\}$ of $G$. Essentially, each clause $(u \vee v)$ is padded with a unique extra variable. 
We prove that if $d-pw(G)$ is at least $k$ then the size of a $d$-SMNBP representing 
$\psi(G)$ is at least $\Omega(c^k)$ where $c=(8/7)^{1/12}$. 
Note that the lower bound is \emph{scalable} in the sense that it does not depend on $d$.
%That is, the lower bound \emph{does not}
%depend on $d$ and, in particular, remains exponential for arbitrary large $d$.  
%In the rest of this section, we refer to it as a \emph{scalable} lower bound. 

We apply the above statement to obtain a lower bound for a concrete class of CNFs: $\psi(K_n)$.
It is easy to see that $\psi(K_n)$ can be represented by a polynomial size $n$-SMNBP
(in fact, by a conjunction of $n$ $1$-NBPs). However, reduction the allowed number of repetitions
from $n$ to $n^{a}$ for an arbitrary constant $0<a<1$ results in an exponential lower bound. 
In particular, combining the $\Omega(c^k)$ lower bound with an upper bound on the number of edges in
a graph of bounded pathwidth, 
we demonstrate that an $n^{a}$-SMNBP representing $\psi(K_n)$ has size at least $\Omega(c^{(n^b)})$.
where $b$ is an arbitrary constant such that $a+b<1$. 

The above lower bound gives rise to several lines of further research.
First, we conjecture that w.r.t. the treewidth alone (without bounded maximal degree assumption), 
$d$-pathwidth behaves just like the ordinary pathwidth: admitting the $\Omega(\log n* tw(G)$ lower bound for
some infinite class of graphs (the constant at the $\Omega$ may depend on $d$).
If this conjecture is confirmed, an XP lower bound in terms treewidth will follow for $d$-SMNBPs for each constant $d$. 
The second research direction is to 'upgrade' the lower bound for $d$-SMNBPs to a lower bound for a more 'mainstream'
model.  
%We find it most intriguing to investigate whether the lower bound holds for non-monotone read $d$-times
%branching programs.  
%To the best of our knowledge, no superpolynomial lower bounds are known for read $n^{\alpha}$ times
%branching programs ($0 \leq \alpha \leq 1$) even if we assume that the branching programs are  Indexed Binary Decision Diagrams
%(IBDDs) \cite{...}.  
As a first step towards this direction, we pose open questions
as to whether the proposed lower bound holds for $d$-MNBPs (without the 'separable' assumption) 
and  whether $\psi(G)$ for $G$ with $2$-pathwidth
at least $k$ requires an exponential in $k$ size representation by a conjunction of two OBDDs.   
The last question aims to investigate the splitting power in the non-monotone case.
The non-monotone splitting power of even read-twice branching programs is greatly enhanced by
the existence of inconsistent paths. This turns the splitting of the set of clauses (which is, essentially,
the case for monotone branching programs) into intersection of the sets of satisfying assignments which can be 
very chaotic and way harder to grasp.  In fact, this is exactly the reason why our argument for the
monotone lower bound fails in the non-monotone case.

\subsection{Motivation}
This line of research is well motivated in several different fields.

Representation of CNFs of bounded treewidth in terms of 'weak' 
classes of Boolean circuits is being actively studied in the area 
of knowledge compilation. A well known fact is the existence of a
'watershed'  between DNNF (Decomposable Negation Normal Forms) based
models and those based on read-once branching programs in the sense that
the former have FPT sized representation of CNFs of bounded primal treewidth
while the latter do not in general. So, it is interesting to see if the 'border' can be crossed
by branching programs through a slightly increased repetition. In this paper we show
that the answer if positive for CNFs with bounded number of occurrences of each variable
but in general the question remains open.

From the graph theoretical perspective, the new parameters, the result 
for bounded degree graphs, and several open questions
contribute to the area of graph sparsity as per the landmark book \cite{sparsitybook}.
The books discusses several variants of vertex and edge colouring.
The graph $c$-pathwidth $k$ are nothing else than graphs whose edge can be coloured
in $c$ colours so that the pathwdith of each monochromatic component  
is at most $k$. 

From the perspective of circuit
complexity, the lower bound proposed in this paper can be seen as a lower bound for a representation of a 'simple'
circuit by a restricted model (this point of view is stated in \cite{WegBook} in the context of FBDDs). 
In case of the concrete lower bound for $\psi(K_n)$, the circuit is so simple that it can be represented by the
very same restricted model with a slightly increased repetition.
From the perspective of future research, the most intriguing question is whether the scalable lower bound
for $\psi(K_n)$ holds without the monotonicity assumption. The positive answer will mean a significant breakthrough
breaking the polylogarithmic repetition barrier for branching programs. On the other hand, the negative answer 
will also mean a significant insight as to how non-monotonicity 'beats' monotonicity in this particular context.

\subsection{Related work}
An FPT upper bound for DNNFs parameterized by treewidth has been obtained in \cite{DarwicheJACM}.
Subsequent research resulted in refinement of the upper bound to several restricted
DNNF classes such as Decision DNNFs \cite{DesDNNF}. On the side of syntactic 
read-once branching programs, $1$-NBP requires, in general, XP size
parameterized by CNFs primal treewidth \cite{RazgonAlgo}.
The result of \cite{RazgonAlgo} also holds for $c$-OBDDs for a constant $c$ \cite{obddtcompsys} (even if they are generalized
to being non-deterministic and semantic). However, the querying pattern of $c$-OBDDs is quite restrictive
and allows establishing lower bounds that (to the best of our knowledge) 
are not currently known for more general models of branching
programs. 

For (syntactic) $k$-NBPs exponential lower bounds are known for $k=O(\log n)$ \cite{readktimes}.
With the extra assumption that the branching programs are deterministic and oblivious,
$k$ can be increased to $O(\log^2 n)$ \cite{logn2}. We are not aware of exponential lower bounds for
higher values of $k$ even if the oblivious branching program is further restricted to $k$-IBDDD \cite{kibdd}.
In fact, we are not aware for such lower bounds even if we restrict a $k$-IBDD to be just
a conjunction of $k$ OBDDs. On the other hand, for $k$-OBDDs, exponential lower bounds are
known for $k=o(n/log n)$ (Corollary 7.5.10 of \cite{WegBook}).

Our definition of monotone NBPs (disallowing negative literals as labels on the edges)
is as in \cite{moncomp} and \cite{RazFCT}. 
It is known that there are problems in monotone NP require monotone circuits 
of exponential size \cite{AlonMono}  and there are problems in monotone P requiring
monotone formulas of exponential size \cite{RobereMono}. The latter lower bound implies
an exponential lower bound for monotone switching and rectifier networks, that is monotone
NBPs with unbounded repetition \cite{RazFCT}.

Several new structural graph parameters have been introduced in the recent years, e.g. \cite{twinwidth1}
and \cite{simwidthpaper}. The parameters introduced in this paper can be seen as a contribution to this trend. 
%An alternative definition of monotone branching program 
%where the monotonicity is defined w.r.t. an ordering of states is considered e.g. in \cite{...}.
%The branching programs as per the latter definition can compute non-monotone functions.
%It remains to consider graph theoretical related works. 
%We need to consider works similar to the Woods' bound and also examples of improper edge
%coloring 

\subsection{Structure of the paper}
Section \ref{sec:prelim} introduces the necessary background.
In Section \ref{sec:graph} we define the new graph parameters
and upper bound them by a function of treewidth and maxdegree
of the considered graph. Next, in Section \ref{sec:obdd}, we use the upper
bound to establish an upper bound for a conjunction of two OBDDs.
Finally, in Section \ref{sec:mono}, we prove a lower bounded for the separable
monotone branching programs of bounded repetition.
Proofs of some statements of Section \ref{sec:mono} are postponed to the
appendix.

\section{Preliminaries} \label{sec:prelim}
\subsection{Models of Boolean functions}
A literal is a Boolean variable or its negation. 
In this paper when we use a set of literals, we mean a \emph{proper}
set of literals where a variable cannot occur along with its negation.
If a variable $x$ occurs in a set $S$ of literals, it can occur \emph{positively},
if $x \in S$ or \emph{negatively}, if $\neg x \in S$. 
We denote by $Var(S)$ the set of variables occurring in $S$.
We also call $S$ an \emph{assignment} to $Var(S)$. 

A \emph{CNF} $\varphi$ is a set of \emph{clauses} and each clause
is just a set of literals.  We denote by $Var(\varphi)$ the set of all
variables occurring in the clauses. A set $S$ of literals \emph{satisfies} a clause $C$
is $S \cap C \neq \emptyset$. $S$ satisfies a set of clauses if it satisfies
each clause of the set. A \emph{satisfying assignment} of a CNF $\varphi$
is an assignment to $Var(\varphi)$ satisfying $\varphi$.

For a Boolean function $F$, we denote by $Var(F)$ the set of variables
of this function. An assignment $S$ of $Var(F)$ is a satisfying assignment
for $F$ if $F$ is $true$ on the tuple where each variable occurring positively
in $S$ is assigned with $true$ and each variable occurring negatively in $S$
is assigned with $false$. 

\begin{definition} [{\bf NBPs}]
A \emph{non-deterministic branching program} (NBP) $Z$
is a directed acyclic graph ($DAG$), multiple edges allowed, with one source and one sink 
some edges of which are labelled with literals. 
We denote by $Var(Z)$ the set of all variables whose literals label the edges of $Z$.  

A (directed) path $P$ of $Z$ is consistent if its labels do not include
a variable along with its negation. For a consistent path $P$, we denote
by $A(P)$ the set of literals labelling the edges of $P$. 

A \emph{satisfying assignment} of $Z$ is a set $S$ of literals with 
$Var(S)=Var(Z)$ and such that there is a consistent source-sink path $P$
with $A(P) \subseteq S$.  In this case, we sometimes say that $P$ \emph{carries} $S$.

We say that $Z$  \emph{represents} a CNF $\varphi$ (respectively, a Boolean function $F$)
if the set of satisfying assignments of $Z$ is the same as that of $\varphi$
(respectively, of $F$).    

We denote by $|Z|$ the number of edges of $Z$. 
\end{definition}

\begin{remark}
There is no point to have two unlabelled edges between the given pair of vertices
$x$ and $y$ or two edges labelled with the same literal. Therefore, we can
assume that the number of multiple edges between the given pair of vertices is at 
most $2n+1$, the total number of literals plus possibly one unlabelled edge.
That is, $|Z|$ and $|V(Z)|$ (the number of vertices of $Z$) are polynomially related.
Therefore,  for the purposes of this paper, we can use either measure.
However, for the upper bound in Section \ref{sec:obdd}, we use $|V(Z)|$ simply
because the bound is based on existing upper bound that also uses $|V(Z)|$.
On the other hand,  the lower bound in Section \ref{sec:mono}
is stated for $|Z|$. The reason is, again, a pure convenience: the lower bound is
proved for an auxiliary branching program having at most $3$ times more edges 
than the original one, while the number or vertices can grow quadratically because
of subdivision of edges. Therefore, the use of $|Z|$ preserves the asymptotical lower bound
for the original branching program.  
\end{remark}

\begin{definition} [{\bf $d$-NBPs}]
A \emph{syntactic read-$d$-times} NBP ($d$-NBP)
is an NBP $Z$ where on each path $P$ and each variable $x \in Var(Z)$,
the number of occurrences of $x$ as a label of an edge of $P$ is at most $d$.  
\end{definition}

%Now, we define two further restrictions of $k$-NBPs that we will use
%in this paper. 
\begin{definition} [{\bf Monotone and separable branching programs}]
An NBP $Z$ is \emph{monotone} if negative literals do not occur as labels
of the edges of $Z$. 
An $d$-NBP $Z$ is \emph{separable} if each source-sink path $P$ can be partitioned
into at most $d$ edge-disjoint \emph{read-once} subpaths.
For example, if the sequence of variables queried along a path is $x_1,x_3,x_1,x_2$
then this path can be partitioned into two read-once subpath whose edges query
variables $x_1,x_3$ and $x_1,x_2$, respectively. However, if the sequence is $x_1,x_3,x_1,x_2,x_2$
then such a partition is not possible. 

We abbreviate a monotone $k$-NBP as $k$-MNBP and a monotone separable $k$-NBP
as $k$-SMNBP.
\end{definition} 

\begin{remark}
The querying constraint imposed by separable $d$-NBPs
is weaker than that of $d$-IBDD \cite{kibdd}.
In particular, separable $d$-NBPs do not place any constraints
on specific orders of querying variables within read-once fragments. 
\end{remark}

\begin{definition} [{\bf OBDDs}]
An ordered binary decision diagram $Z$ is a DAG, multiple edges allowed, with a single source
and two sinks, one labelled with $True$, the other labelled with $False$.
Each non-sink node has two outgoing edges labelled with opposite literals
of the same variable. The labelling of the edges is \emph{read-once}:
for each path $P$ of $Z$ there are no two different edges labelled by literals of the
same variable. The labelling is also \emph{oblivious}: there is a permutation 
$\pi$ of variables ( treated as a linear order) : for each path $P$ whenever a literal
of $y$ occurs on $P$ after a literal of $x$, it holds that $y$ occurs after $x$ in $\pi$.

We denote by $Var(Z)$ the set of variables whose literals label the edges of $Z$. 
We denote by $A(P)$  the set of literals labelling the edges of a path $P$ of $Z$.
The function $F(Z)$ represented by $Z$  is a function whose set of variables is $Var(Z)$
and the set of satisfying assignments consists of precisely those $S$ such that there is
a path $P$ from the source to the $True$ sink such that $A(P) \subseteq S$.
\end{definition}

\begin{definition}
Let $Z_1, \dots, Z_q$ be OBDDs.
The function $F=F(Z_1) \wedge \dots \wedge F(Z_q)$
is called the \emph{conjunction} of $Z_1, \dots, Z_q$ 
\end{definition}

\begin{remark}
The conjunction of OBDDs $Z_1, \dots, Z_q$ can be easily represented
as $q$-NBP as follows. Transform each $Z_i$ into a $1$-NBP $Z'_i$ by removal of the $False$
sink and all the nodes from which the $True$ sink cannot be reached. 
Then for each $1 \leq i \leq q-1$ identify the sink of $Z_i$ with the source of $Z_{i+1}$.
%See Figure \ref{...} for the illustration. 
The same 'chaining' approach but with a slightly
more tedious implementation can be used to demonstrate that the conjunction
of $q$ OBDDs can be represented as a $q$-IBDD, a restricted class of deterministic
read $q$ times branching programs.  
\end{remark}

\subsection{Graphs and their structural parameters}
We use a standard terminology related to graphs as in e.g. \cite{Diestel3}.
In particular, we denote by $G[S]$ the subgraph of $G$ induced by $S \subseteq V(G)$.

\begin{definition} [{\bf Treewidth and pathwidth}]
A \emph{tree decomposition} of a graph $G$ is a pair $(T,{\bf B})$ where $T$ is a
tree and ${\bf B}$ is a set of bags $B_x$ corresponding to the nodes $x$ of $T$.
Each bag is a subset of $V(G)$ and the following conditions must be met:
(i) \emph{union}, that is $\bigcup_{x \in V(T)} B_x=V(G)$, (ii) \emph{containment},
that is for each $e \in E(G)$ there is $x \in V(T)$ such that $e \subseteq B_x$, and
(iii) \emph{connectedness}, that is for each $u \in V(G)$, the set $\{x| u \in B_x\}$ induces
a connected subgraph of $T$. 

If $T$ is a path then $(T,{\bf B})$ is called a \emph{path decomposition} of $G$. 

The width of $(T,{\bf B})$ is the size of the largest bag minus one.
The \emph{treewidth} of $G$, denoted by $tw(G)$ is the smallest width of a tree decomposition of $G$.
The \emph{pathwidth} of $G$, denoted by $pw(G)$ is the smallest width of a path decomposition of $G$. 
\end{definition} 

We conclude this section with three facts about treewidth and pathwidth along with
literature references for relevant proofs.

\begin{proposition}  \label{prop:clique}
If $S$ is a clique of $G$ then $S$ s a subset of a bag in every tree decomposition of $G$. 
\end{proposition}

\begin{proposition} \label{prop:nk}
A graph of treewidth $k$ has at most $nk$ edges.
\end{proposition}

\begin{proposition} \label{prop:log}
There is an infinite class $\mathcal{G}$ of graphs for which there is a constant $\alpha$
such that for each $G \in \mathcal{G}$, $pw(G) \geq \alpha*tw(G)*\log n$.
\end{proposition}

Proposition \ref{prop:clique} appears in \cite{Diestel3}  as Lemma 12.3.5,
Proposition \ref{prop:nk} appears in \cite{Reedtutorial} as statement 1.10, 
A class as stated in Proposition \ref{prop:log} with an additional property
that the max-degree of all graphs is $5$ is provided in \cite{RazgonAlgo}.

\section{New parameters and their upper bound for graphs of bounded degree} \label{sec:graph}

%The definition below formalizes the notion of edge colouring
%so that each monochromatic component has a low pathwidth. 

\begin{definition}
Let $d \geq 1$ be an integer and $G$ be a graph.
The $d$-\emph{pathwidth} of $G$ denoted by $d-pw(G)$is the smallest $k$ such that
there are subgraphs $G_1, \dots G_d$ of $G$ each of pathwidth
at most $k$ and such that $G=G_1 \cup \dots \cup G_d$. 

The \emph{clique preserving} $d$-pathwidth denoted by $d-cpw(G)$
is defined analogously with the only extra requirement that each
complete subgraph of $G$ is a subgrpah of some $G_i$.
\end{definition}

\begin{example}
A rectangular grid has $2$-pathwidth $1$. Indeed, let one subgraph be induced by all the 'horizontal'
edges and the other subgraph be induced by all the 'vertical' edges. This way the grid is represented
as the union of two subgraphs each connected component of each subgraph is a path.
\end{example}

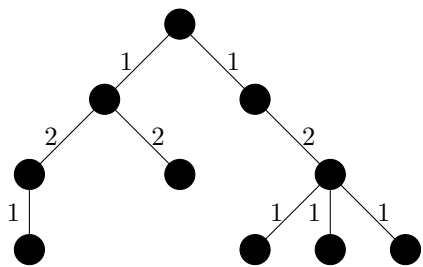
\begin{figure}[h]
\begin{tikzpicture}
\draw [fill=black]  (4,4) circle [radius=0.2];
\draw [fill=black]  (3,3) circle [radius=0.2];
\draw [fill=black]  (5,3) circle [radius=0.2];
\draw [fill=black]  (2,2) circle [radius=0.2];
\draw [fill=black]  (4,2) circle [radius=0.2];
\draw [fill=black]  (6,2) circle [radius=0.2];
\draw [fill=black]  (2,1) circle [radius=0.2];
\draw [fill=black]  (5,1) circle [radius=0.2];
\draw [fill=black]  (6,1) circle [radius=0.2];
\draw [fill=black]  (7,1) circle [radius=0.2];

\draw (4,4) --(3,3); 
\draw (4,4) --(5,3);
\draw (3,3) --(2,2);
\draw (3,3) --(4,2);
\draw (5,3) --(6,2);
\draw (2,2) --(2,1);
\draw (6,2) --(5,1);
\draw (6,2) --(6,1);
\draw (6,2) --(7,1);
\node [left]  at (3.5,3.5)  {$1$};
\node [right]  at (4.5,3.5)  {$1$};
\node [left]  at (2.5,2.5)  {$2$};
\node [right]  at (3.5,2.5)  {$2$};
\node [right]  at (5.5,2.5)  {$2$};
\node [left]  at (2,1.5)  {$1$};
\node [left]  at (5.5,1.5)  {$1$};
\node [left]  at (6,1.5)  {$1$}; 
\node [right]  at (6.5,1.5)  {$1$};       
\end{tikzpicture}
\caption{A tree as the union of two subgraphs of pathwidth $1$. The edges of one subgraph are labelled
with $1$, the edges of the other subgraph are labelled with $2$. }
\label{treegraph}
\end{figure}

\begin{example} \label{extree}
A tree has $2$-pathwidth $1$.
Indeed, let $T$ be a tree. Pick an arbitrary node $r$ of $T$ and let it be the root.
Then the edges are naturally divided into layers. The edges between the root and its children are
of layer one. the edges between the children of the root and their children and layer $2$ and so on. 
Let $G_1$ and $G_2$ be the subgraphs of $G$ induced by the edges of the odd and even layers, respectively. 
Then each connected component of each $G_i$ is a star and hence both subgraph have pathwidth $1$.
This approach is demonstrated in Figure \ref{treegraph}.
\end{example}

The main question studied in this section is the following:
can $d-pw(G)$ and $d-cpw(G)$ be upper bounded by a function
of $tw(G)$, the treewidth of $G$ for any constant $d$? 
For graphs of a bounded degree the answer is positive already for $d=2$,
for graphs in general we will conjecture that this is not the case in a strong
sense.

In order to address the case of a bounded degree we need
the notion of a tree partition width.

\begin{definition}
Tree-partition decomposition of a graph $G$ is a pair $(T,{\bf B})$ where
$T$ is a forest and ${\bf B}$ is a set of bags corresponding to the nodes of $T$
that constitute a \emph{partition} of $V(G)$.  Furthermore, let $t_1,t_2$ be two
distinct nodes of $T$. Then there $B(t_1)$ and $B(t_2)$ are adjacent in $G$
if and only if $t_1$ and $t_2$are adjacent in $T$. 
The width of $(T,{\bf B})$ is the largest size of a bag.
The tree-partition width of $G$ denoted   
by $tpw(G)$ is the smallest width of a tree-partition decomposition of $G$. 
\end{definition}

The tree partition width of $G$ can be linearly upper bounded
by the product of $tw(G)$ and $\Delta(G)$ \cite{treepartition}.  %Wood 'On tree-partition-width'.

\begin{theorem} \label{mainbound}
There is a constant $\gamma$ such that for each graph $G$ with at least one edge, $tpw(G) \leq \gamma*\Delta(G)*tw(G)$,
where $\Delta(G)$ is the max-degree of $G$.  
\end{theorem}

We are going to show that $2-cpw(G) \leq 2*tpw(G)-1$.
Combined with Theorem \ref{mainbound}, this implies 
that $2-cpw(G)=O(\Delta(G)(tw(G))$.
In order to upper bound $2-cpw(G)$ in terms of $tpw(G)$
we turn the witnessing tree-partition decomposition for
$tpw(G)$ into a tree decomposition by leaving the same underlying tree
with arbitrarily identified root and by adding the vertices of the bag of the parent to each non-root bag.
The nodes of the underlying rooted tree are naturally defined into layers.
We enumerate the layers from the top to the bottom.
Then we create two graphs one induced by the union of bags 
of even layers and one induced the union of bags at odd layers.
We then show that these two graphs witness the bounded $2-cpw(G)$.
Example \ref{extree} and Figure \ref{treegraph} illustrate this approach in case $G$ is a tree. 
A formal description is provided in the theorem below.

\begin{theorem}  \label{mainupperbound}
$2-cpw(G) \leq 2tpw(G)-1$.
In particular, it follows from Theorem \ref{mainbound}
that $2-cpw(G) \leq 2\gamma(\Delta(G)*tw(G)$.
\end{theorem}

{\bf Proof.}
Let $(T,{\bf B})$ be a tree-partition decomposition of width $tpw(G)$.
Turn $T$ into a rooted tree by arbitrarily picking a node of $T$
and calling it $root(T)$.  For each non-root node $t \in V(T)$,
let $p(t)$ be the parent of $t$. 
Let $(T,{\bf B^*})$ be a pair where ${\bf B^*}$  is a set of bags 
associated with the nodes of $T$ so that $B^*(root(T))=B(root(T))$
and, for every other $t \in V(T)$, $B^*(t)=B(t) \cup B(p(t))$.
It is not hard to see that $(T,{\bf B^*})$ is a tree decomposition of $G$. 
We call it the tree decomposition \emph{induced} by $(T,{\bf B})$.

Next, we partition the nodes of $T$ into \emph{layers} as follows:
$root(T)$ is the only node of layer $1$, the children of the root are the nodes
of layer $2$, the children of children are the nodes of layer $3$ and so on.
Put it differently, the layer number of a node is its distance from the root plus one. 
We denote by $Even(T)$ and $Odd(T)$ the sets of nodes that belong to the even
and odd layers respectively. 
Let $V_{even}=\bigcup_{t \in Event(T)} B^*(t)$ and 
$V_{odd}=\bigcup_{t \in Odd(T)} B^*(t)$.
We denote $G[V_{even}]$ and $G[V_{odd}]$ by $G_{even}$ and $G_{odd}$ and call
them the \emph{even} and \emph{odd} subgraphs of $G$. 

We claim that $G_{even}$ and $G_{odd}$ are the subgraphs witnessing
$2-cpw(G) \leq 2tpw(G)-1$. Indeed, 
by Proposition \ref{prop:clique}, each complete subgraph of $G$ is a subgraph of some $G[B^*(t)]$
and hence, by construction, a subgraph of either $G_{even}$ or $G_{odd}$.
It remains to show that the pathwidth of $G_{even}$ and $G_{odd}$ is at 
most $2tpw(G)-1$.  We demonstrate this for $G_{odd}$, the proof for $G_{even}$
is symmetric. 

We construct a path decomposition of $G_{odd}$ whose bags
are $B^*(t)$ for each $t \in Odd(t)$ 
By construction, the size of each bag is at most $2tpw(G)$. 
It remains to be shown that the bags can be linearly 
ordered so that the bags containining each vertex of $G$
form an interval. 

Partition $Odd(T)$ into subsets $\{S_1, \dots, S_q\}$
such that two nodes get into the same subset if and 
only if they are siblings in $T$. Arbitrary order each $S_i$ into
a sequence and denote it
by $\pi_i$. Let $\pi=\pi_1+ \dots+\pi_q$.
With a slight abuse of notation, we treat $\pi$ as a path with an edge between
every pair of consecutive elements. 
We claim that  $(\pi,{\bf B_0})$ where for each $t \in Odd(t)$,
$B_0(t)=B^*(t)$ is a path decomposition of $G_{odd}$. As said in the previous
paragraph, it remains to verify the connectedness property.

Let $u \in V(G_{odd})$.  Let $t \in T$ such that $u \in B(t)$.
If $t \in Odd(T)$ then in $(T,{\bf B^*})$ $u$ belongs to $B^*(t)$ and,
possibly, to the bags of the children of $t$ that all belong to $Even(t)$.
It follows that the only bag of ${\bf B_0}$ containing $u$ is $B_0(t)$
and hence the connectedness clearly holds. 
If $t \in Even(T)$ then in $(T,{\bf B^*})$, apart from $B^*(t)$,
$u$ also belongs to the bags of the children of $t$ that are
all siblings and form some $S_i$. The bags of these $S_i$ in ${\bf B_0}$ are 
the only bags containing $u$. By construction they form an interval.
$\blacksquare$

%Check in the paper of Wood whether $T$ can be assumed binary
%introduce the related discussion.
As $c$-pathwidth never exceeds its clique preserving variant,
the above upper bound holds for $2-pw(G)$ as well. 

We do not know whether for a constant $c$, $c-pw(G)$ can be upper-bounded
by a function of $tw(G)$ alone. Moreover, we are not aware of existing
results supporting intuition that this might be the case.
On the other hand, it is known that $tpw(G)$ cannot be upper bounded 
by a function of $tw(G)$ alone \cite{treepartition}. We therefore conjecture that the lower bound
as in Proposition \ref{prop:log} also holds for $c$-pathwidth.

\begin{conjecture} \label{largepw}
For each integer $c \geq 1$ there is a constant $\alpha_c$ and a 
class ${\bf G}_c$ of graphs of unbounded treewidth such that
for each $G \in {\bf G}_c$, $c-pw(G) \geq \alpha_c*tw(G)*\log n$.
\end{conjecture}

\begin{remark}
It is not hard to see that the connected components of both graphs $G_{even}$
and $G_{odd}$, as in the proof of Theorem \ref{mainupperbound},
are subgraphs of $G$ induced by $B(t)$ for some $t \in V(T)$ and the respective bags
of the children of $t$. Let $P$ be a path of such a component.
As in $(T, {\bf B})$, two vertices in the bags of distinct children of $t$ are not adjacent (by
the properties of the tree-partition width)
and the path of vertices of the same bag is of size at most $tpw(G)$,
any subpath of $P$ of length $tpw(G)+1$ contains a vertex of $B(t)$.
Therefore the length of $P$ is at most $(tpw(G)+1)*tpw(G) \leq O(\bigtriangleup^2(G)*tw^2(G))$
by Theorem \ref{mainbound}. In other words the path length in each $G_{even}$ and $G_{odd}$
is upper bounded by a function of the max-degree and the treewidth of $G$.
Therefore, a reasonable first step towards resolving Conjecture \ref{largepw} would be to design
a class of graphs in which there is no two colouring of edges with the length of monochromatic 
paths upper bounded by a function of the treewidth. 
\end{remark}

\section{Bounded treewiwdth and degree CNFs and conjunction of OBDDs}  \label{sec:obdd} 
Throughout this section $\varphi$ is a CNF and $|Var(\varphi)|$ is denoted by $n$. 
\begin{definition}
The \emph{primal}  graph $G_{\varphi}$ has $Var(\varphi)$ as the set of vertices.
Two variables are adjacent in $G_{\varphi}$ if and only if they occur in the same
clause of $\varphi$. 
The primal treewidth and pathwidth of $\varphi$ are respective treewidth
and pathwdith of $G_{\varphi}$ and are denoted by $tw(\varphi)$ and $pw(\varphi)$,
respectively.  
\end{definition}

By analogy with graphs, we introduce the notation
$\Delta(\varphi)=max_{x \in Var(\varphi)} d_{\varphi}(x)$
where $d_{\varphi}(x)$ is
the number of clauses of $\varphi$ where the variable $x$ occurs. 
In this section we show $\varphi$ can be represented as a conjunction
of two OBDDs of size at most $2^{O(\Delta(\varphi)*tw(\varphi))}*n$.  
For the proof we will use Theorem \ref{mainupperbound}  (in particular, we will
clarify why we need the clique preserving variant of $2$-pathwidth) and the following
result from \cite{VardiTWD}.

\begin{theorem} \label{pwupper}
A CNF $\psi$ can be represented by an OBDD of at most $2^{pw(\psi)}*|Var(\psi)|$ nodes.
\end{theorem}

\begin{theorem}  \label{twoobddupper}
A CNF $\varphi$ can be represented
as the conjunction of two OBDDs of at most $2^{2\gamma*\Delta(\varphi)*tw(\varphi)^2}*n$ nodes
each where $\gamma$ is the constant as in Theorem \ref{mainupperbound}.
\end{theorem}

{\bf Proof.}
We demonstrate the existence of two CNFs 
$\varphi_1$ and $\varphi_2$ such that $\varphi_1 \wedge \varphi_2=\varphi$ and 
for each $i \in {1,2}$, $pw(\varphi_i) \leq 2\gamma*\Delta(\varphi)*tw(\varphi)^2$
%the pathwidth of the primal graph of
%each $\varphi_1$ is at most 
%$2\gamma*d*k^2$. 
Then we apply Theorem \ref{pwupper}.

Recall that $G_{\varphi}$ denotes the primal graph of $\varphi$. 
As in the primal graph clauses turn into cliques, and each clique, in turn, is a subset of
some bag by Proposition \ref{prop:clique},  no clause can be of size larger than $tw(\varphi)+1$.
Thus each variable of $\varphi$ can be in the same clause with at most $\Delta(\varphi)*tw(\varphi)$
other variables and hence this an upper bound on the degree of $G_{\varphi}$.

By Theorem \ref{mainupperbound}, there are two graphs $G_1$ and $G_2$
such that $G_1 \cup G_2=G_{\varphi}$, each clique of $G$ is a subgraph of some $G_i$
and the pathwidth of each $G_i$ is at most $2\gamma*\Delta(\varphi)*tw(\varphi)^2$.

Create CNFs $\varphi_1$ and $\varphi_2$ as follows.
For each clause $C$ of $\varphi$, if the clique induced by $C$
is a subgraph of $G_1$, let $C$ be a clause of $\varphi_1$
otherwise %(when, by the clique preservation $G[C]$ is a subgraph of $G_2$!)
let $C$ be a clause of $\varphi_2$. 

By construction,  $\varphi_1 \wedge \varphi_2=\varphi$ and 
$G_{\varphi_1}$ is a subgraph of $G_1$. 
As for each clause $C$, $G[C]$ is a clique, if $G[C]$ is not a subgraph of $G_1$
then, by clique preservation \footnote{this is why the clique preserving variant
of $2$-pathwidth is needed!}, $G[C]$ is a subgraph of $G_2$.
Therefore $G_{\varphi_2}$ is a subgraph of $G_2$.
We conclude that both $pw(\varphi_1)$ and $pw(\varphi_2)$
are at most $2\gamma*\Delta(\varphi)*tw(\varphi)^2$.
$\blacksquare$

\section{Lower bounds depending on $d$-pathwidth} \label{sec:mono}
Recall from the Preliminaries section that a $d$-MNBP $Z$
is \emph{separable}, abbreviated as $d$-SMNBP, if every source-sink path of $Z$ can be
partitioned into at most $d$ read-once subpaths.

In this section we prove a lower bound for $d$-SMNBPs representing
a class of CNFs. 
%In this section we define a quite general subclass of
%\emph{monotone} read $c$-times non-deterministic branching 
%programs ($c$-MNBP) and prove a lower bound on the size of these
%branching programs representing a class of CNFs. 
The lower bound is exponential is terms of $d$-pathwidth of 
the primal graphs of these CNFs. 

%First of all, we need to define the considered class of branching programs.
%\begin{definition}
%We call $Z$ a \emph{separated} $c$-MNBP ($c$-SMNBP) if 
%$Z$ is a $c$-MNBP in which every source-sink path can be partitioned 
%into at most $c$ read-once paths.
%\end{definition}

The considered class of CNFs are in a bijective correspondence
with graphs without isolated vertices. In particular, the 
CNF $\psi(G)$ corresponding to a graph $G$ has $V(G) \cup E(G)$
as the set of variables. 
The variables of $V(G)$ and $E(G)$ are, respectively, the \emph{vertex} and \emph{edge} variables. 
The clauses correspond to $E(G)$.
In particular for each edge $e=\{u,v\}$ of $G$, the corresponding
clause is $(u \vee e \vee v)$.
Note that each $e$ occurs as a variable only in the clause corresponding
to $e$. The role of the edge variables is \emph{padding} that allows any
assignment to the vertex variables to be extended to a satisfying
assignment of $\psi(G)$.
The primal graph $H$ of $\psi(G)$ is obtained from
$G$ by introducing an individual vertex for each $e=\{u,v\}$ 
and making it adjacent to vertices $u$ and $v$. 

\begin{proposition} \label{getridofh}
$d-pw(G) \leq d-pw(H) \leq d-pw(H)+1$.
%It is not hard to see that the $c$-pathwidth of $H$ 
%is at most the $c$-pathwdith of $G$ plus one. 
%Therefore, an exponential lower bound in terms of  
%of the $c$-pathwidth of $G$ implies an exponential lower
%bound in terms of $c$-pathwidth of $H$.
\end{proposition}

{\bf Proof.}
As $G$ is a subgraph of $H$, the first inequality is immediate.
For the other inequality, let $G_1, \dots, G_d$ be subgraphs of $G$ of pathwidth
at most $d-pw(G)$ each whose union is $G$. 
We transform each $G_i$ into $H_i$ as follows. For each $\{u,v\} \in E(G_i)$,
add the unique new vertex $x$ whose neighbours are $u$ and $v$ along with the edges
connecting $x$ to $u$ and $v$.  As each $\{u,v\} \in E(G)$ is an edge of some $G_i$,
the new vertex $x$ whose neighbours are $u$ and $v$ belongs to $V(H_i)$ and the two
edges adjacent to $x$ belong to $E(H_i)$. We conclude that $H_1 \cup \dots \cup H_d=H$.
It remains to show that for each $i \in \{1, \dots, d\}$, $pw(H_i) \leq pw(G_i)+1$.

Let $(P,{\bf B})$ be a path decomposition of $G_i$ having the smallest possible width.
For each $x \in V(P)$ let $r(x)=|E(G_i[B_x])|$.
Put it differently, $r(x)$ is the number of edges of $G_i$ between vertices of $B_x$.
Form a new path $P'$ by replacing each $x$ with a sequence of $r(x)$ nodes and let $B(x)$ be the bag of each node.
Thus the bags $B(x)$ are now in a bijective correspondence with $E(G_i[B_x])$. 
Now, add to each bag $B(x)$ the new vertex of $H_i$ corresponding to the edge of $E(G_i[B_x])$
that corresponds to this bag. This way all the vertices of $H_i$ are accommodated an each
bag becomes larger by at most one element. Each new vertex belongs to exactly one bag
and for each old vertex $u$ the subpath of the nodes of $P$ whose bags contain $u$ may become
longer in $P'$  but is still a subpath. Hence, we have obtained a path
decomposition of $H_i$ of width at most $pw(G_i)+1$.
$\blacksquare$

%Insert the proof
In light of Proposition \ref{getridofh}, we use $d-pw(G)$ rather than $d-pw(H)$ in
the main theorem of this section and the lower bound in terms of $d-pw(H)$ readily follows.  

\begin{theorem}  \label{mainlower}
%There is a constant $\beta>1$ such
%that the following is true.  
Let $d,k \geq 1$ be integers. 
Let $G$ be a graph with $d-pw(G) \geq k$. 
Let $Z$ be a $d$-SMNBP with representing $\psi(G)$.
Then %$|Z| \geq [(8/7)^{1/3}]^{\lfloor k/4 \rfloor}$ or, to put it differently, 
$|Z| \geq \Omega(\beta^k)$ where $\beta=(8/7)^{1/12}$. 
In particular, subject to Conjecture \ref{largepw} being true,
for every constant $d$,  
there is no FPT-sized $d$-SMNBP representation of CNFs of bounded treewidth.
\end{theorem}

%Here we need to provide a proof that 
%for $\psi(K_n)$ a large repetition SMNBP will still have an exonential size. 

%The other implication is an unconditional lower bound on SMNBPs
%with a sublinear repetition represening $\psi(K_n)$.
%For any posiive $\alpha, \beta$ such that $\alpha+2\beta<1$,
%a simple counting argument shows that the $n^{\alpha}$ -pathwdith
%of $K_n$ is at least $n^{\beta}$. Therefore, for a sufficiently large $n$, 
%the size of $n^{\alpha}$-SMNBP is at least $c^{n^{\beta}}$.

Before proving Theorem \ref{mainlower},
we demonstrate its application by proving a lower bound for CNFs $\psi(K_n)$
that is tight in the sense described below.
 
\begin{theorem}
Let $a$ and $b$ be positive constants such that $a+b<1$. 
Then, for a sufficiently large $n$,  $n^a$-SMNBP representing $K_n$ has size at least $\Omega(\beta^{n^{b}})$
where $\beta$  is as in Theorem \ref{mainlower}.

On the other hand,  $\psi(K_n)$ has an $O(n^2)$ representation as $n$-SMNBP
\end{theorem}

{\bf Proof.}
Represent $K_n$ as the union of graphs $G_1, \dots, G_q$ so that $q \leq n^a$.
Then at least one $G_i$ will have at least ${n \choose 2}/n^a$ edges.
For a sufficiently large $n$,  ${n \choose 2}/n^a>n*n^b$.
Then, by Proposition \ref{prop:nk}, the pathwidth of $G_i$ is greater than $n^b$.
We conclude that the $n^a-pw(K_n)>n^b$. 
The  lower bound as specified in the statement immediately follows from Theorem \ref{mainlower}.

For the upper bound, represent $K_n$ as the union of $n$  stars $K_{1,n-1}$,
represent each copy of $K_{1,n-1}$ as $1$-MNBP, the resulting $n$-SMNBP is just their conjunction. 
It remains to show
how to represent $K_{1,n-1}$ as a $1$-MNBP of size $O(n)$.

Let $v_0, \dots v_{n-1}$ be the vertices of $K_{1,n-1}$, $v_0$ being the centre.
Let $e_1, \dots, e_{n-1}$ be the edges connecting $v_0$ to $v_1, \dots, v_{n-1}$ respectively.
That is $\psi(K_{1,n-1})=(v_0 \vee e_1 \vee v_1) \wedge \dots \wedge (v_0 \vee v_{n-1})$

Let $Z'$ be an $1$-MNBP with vertices $x_0, \dots, x_{n-1}$ with $x_0$ being the source and $x_{n-1}$ being
the sink. Introduce an edge from $x_0$ to $x_{n-1}$ and label it with $v_0$.
Then  for each $1 \leq i \leq n-1$ introduce a pair of parallel edges between $x_{i-1}$ and $x_i$
label one of them with $v_i$ and the other with $e_i$.
A direct inspection shows that $Z'$ represents 
$\blacksquare$

In order to prove Theorem \ref{mainlower},
we introduce a number of auxiliary statements.
Their proofs are provided in the appendix (but one that is provided in this 
section).

First of all, 
we introduce a restricted version of $d$-SMNBP 
called $d$-SMNBP \emph{with yardsticks}.
We then show that the version with yardsticks
simulates the $d$-SMNBP with only a linear increase
in the number of edges. Then we prove Theorem \ref{mainlower}
under assumption that the underlying $d$-SMNBP 
is with yardsticks.  Theorem \ref{mainlower} without
the assumption will immediately follow from the 
combination of these two statements.  

\begin{definition}
Let $Z$ be a $d$-SMNBP.
We say that $Z$ \emph{has yardsticks}
if every path $P$  has $a+1$ different vertices 
$u_a, \dots, u_{a+1}$, $a \leq d$ where $u_1$ is the source of $Z$,
$u_{a+1}$ is the sink of $Z$ such that for every  $1 \leq i \leq a$,
any $u_i,u_{i+1}$ path of $Z$ is read once.
The sets $u_1, \dots, u_{a+1}$ are called the \emph{yardsticks} 
(note that $P$ may have several sets of yardsticks). 
\end{definition}

\begin{figure}[h]
\begin{tikzpicture}
\draw [fill=black]  (5,5) circle [radius=0.2];
\draw [fill=black]  (4,4) circle [radius=0.2];
\draw [fill=black]  (6,4) circle [radius=0.2];
\draw [fill=black]  (5,3) circle [radius=0.2];
\draw [fill=black]  (5,2) circle [radius=0.2];
\draw [fill=black]  (5,1) circle [radius=0.2];

\draw [->] (5,5) --(4.15,4.15);
\draw [->] (5,5) --(5.85,4.15);
\draw [->] (4,4) --(4.85,3.15);
\draw [->] (6,4) --(5.15,3.15);
\draw [->] (5,3) --(5,2.2);  
\draw [->] (5,2) --(5,1.2);

\node [left]  at (4.5,4.5)  {$v_3$};
\node [right]  at (5.5,4.5)  {$v_1$};
\node [left]  at (4.5,3.5)  {$v_3$};
\node [right]  at (5.5,3.5)  {$v_2$};
\node [left]  at (5,2.5)  {$v_1$};
\node [left]  at (5,1.5)  {$v_2$};

\node [above]  at (5,5.1)  {$x_1$};
\node [right]  at (6.1,4)  {$x_2$};
\node [right]  at (5.1,3)  {$x_3$};

\end{tikzpicture}
\caption{A $2$-SMNB without yardsticks}
\label{noyard}
\end{figure}
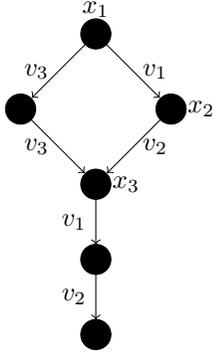

\begin{example}
Consider the $2$-SMNBP on Figure \ref{noyard}.
The variables labelling its edges $v_1,v_2,v_3$.
Also, $x_1, x_2, x_3$ are vertex names we need for further reasoning. 
Let $P$ be the source sink path with edges labelled by $v_1,v_2,v_1,v_2$. 
The only way to partition $P$ into two read-once fragments if to have 
the first fragment consisting of the first two edges and the second fragment
consisting of the last two edges. In other words, the first and the last nodes
of the first fragment are $x_1$ and $x_3$, respectively. 
However, there is another path between $x_1$ and $x_3$ both edges of which are labelled
with $v_3$, that is this alternative path is not read-once. Therefore, path $P$ cannot be assigned
with yardsticks. 

This $2$-SMNBP can be easily turned into one with yardsticks by subdivision of the edge $(x_2,x_3)$.
In particular, introduce a new vertex $x_4$ and replace $(x_2,x_3)$ with two new edges $(x_2,x_4)$ and 
$(x_4,x_3)$. Label $(x_2,x_4)$ with $v_2$.  The resulting branching program represents the same function
as the original one but both source-sink paths have yardsticks.
%: the node between two edges labelled with $v_3$ and
%$x_4$, respectively.

The simulation used for the proof of Theorem \ref{yardsimulate} uses essentially the same approach. 
\end{example}

\begin{theorem} \label{yardsimulate}
A $d$-SMNBP can be simulated by an $d$-SMBP with yardsticks
with at most three times more edges.
\end{theorem}

We proceed to discuss the proof of Theorem \ref{mainlower}
under assumption that $Z$ is an $d$-SMNBP with yardsticks. 
We introduce a probability
space over satisfying assignments of $\psi(G)$
We then prove that the probability of a set of satisfying assignments
satisfying a certain property is at most $(7/8)^{k/4}$.
Next we prove that for each source-sink path 
of $Z$, there are 3 vertices, so that the set of satisfying assignments
carried by the source-sink paths passing through all these three vertices %define carried by as A(P) \subseteq S
satisfies the above property. Combining the $(7/8)^{k/4}$ with the union
bound implies that the number of such triples of vertices is
at least $(8/7)^{k/4}$ meaning that the total number of vertices
(and hence the number of edges) is lower bounded by $\Omega((8/7)^{k/12})$
as required.  

To proceed, let us denote by ${\bf SAT}(G)$ the set of all satisfying assigments
of $\psi(G)$ (recall that, by definition, 
the set of variables of a satisfying assignment of $\psi(G)$ is always $Var(\psi(G))$).

\begin{definition} \label{def:vertedge}
The \emph{Vertex-Edge} probability space of a graph
$G$ denoted by $\mathcal{VE}(G)$
is a probability space whose universe is ${\bf SAT}(G)$.
The probabilities of assignments are defined
as follows.  Let $S \in {\bf SAT}(G)$. 
%Let $SV$ and $SE$ be respective projections of $S$ to the
%sets of vertex and edge variables of $\psi(G)$. 
Let us call an edge $\{u,v\}$  \emph{free} by $S$ if either $u$ or $v$
occur positively in $S$; otherwise, the edge is called \emph{enforced}
by $S$. Let $Free(S)$ be the set of free edges by $S$.
Then $Pr_{\mathcal{VE}(G)}(S)=2^{-(|V(G)|+|Free(S)|)}$ 
($\mathcal{VE}(G)$ is the only probability space we use in this paper. 
Therefore, in what follows we will omit the subscript of $Pr$). 
\end{definition}

It is not hard to observe that $\mathcal{VE}(G)$ is 
indeed a probability space. Indeed, let 
$SV$ be an assignment with $Var(SV)=V(G)$. 
Let ${\bf S}(SV)$ be the set of all $S \in {\bf SAT}(G)$
whose projection to the vertex variables is $SV$.
Then $Free=Free(S)$ is the same for all $S \in {\bf S}(SV)$
(completely determined by $SV$)
and hence $|{\bf S}|=2^{|Free(S)|}$ (the enforced edges assigned positively, the free edges
assigned arbitrarily). 
Then $Pr({\bf S}(SV))=2^{-|V(G)|}$.
As the set of satisfying assignments is the disjoint union of all ${\bf S}(SV)$,
we conclude that the sum of all the probabilities is $1$.

\begin{definition}
Let ${\bf S} \subseteq {\bf SAT}(G)$.
We say that ${\bf S}$ \emph{fixes} a clause $C$ of $\psi(G)$
if there is $C' \subset C$ such that for each $S \in {\bf S}$,
$C' \cap S \neq \emptyset$.  
We can also say that ${\bf S}$ fixes $C$ \emph{with} $C'$ if a specific
subset is needed in the context. 
We say that ${\bf S}$ fixes a set of clauses
if ${\bf S}$ fixes each clause of the set. 
\end{definition} 

\begin{example}
Consider $P_4$, a path of $4$ vertices 
with $v_1,v_2,v_3,v_4$ being the vertices 
and $e_1=\{v_1,v_2\}$, $e_2=\{v_2,v_3\}$, $e_3=\{v_3,v_4\}$
being the edges. 
That is, $\psi(P_4)=(v_1 \vee e_1 \vee v_2) \wedge 
(v_2 \vee e_2 \vee v_3) \wedge (v_3 \vee e_3 \vee v_4)$. 
Let us define a set ${\bf S}$ of satisfying assignments of $\psi(P_4)$ 
as follows.

${\bf S}=\{\{\neg v_1, \neg e_1, v_2, \neg e_2, \neg v_3, e_3, \neg v_4\}, \\
               \{\neg v_1,e_1, \neg v_2, e_2, \neg v_3,e_3, \neg v_4\}, 
                \{\neg v_1,e_1, \neg v_2, \neg e_2,  v_3, \neg e_3, \neg v_4\}         
\}$

Then ${\bf S}$ fixes $(v_1 \vee e_1 \vee v_2)$ with $\{e_1,v_2\}$
and $(v_3 \vee e_3 \vee v_4)$  with $\{v_3,e_3\}$.
However, ${\bf S}$ does not fix $(v_2 \vee e_2 \vee v_3)$
as every proper subset of the clause of falsified by some assignment of ${\bf S}$.   
\end{example}
A \emph{matching} of clauses of $\psi(G)$ is a set of clauses
whose corresponding edges form a matching.  

\begin{theorem} \label{fixedprob}
%There is a constant $0<\alpha<1$
%such that 
For any matching $M$ of
$\psi(G)$ and any set ${\bf S}$ of satisfying assignments of
$\psi(G)$,
that fixes $M$, $Pr({\bf S}) \leq (7/8)^{|M|}$.
%($\mathcal{VE}(G)$ is the only probability space considered in this
%theorem so we do not explicitly specify itas a subscript of $Pr$.)
\end{theorem}

Now we are going to show that for any source-sink path 
of $Z$,  there are three vertices, so that the set of satisfying assignments
carried out by paths going through all these three vertices
fixes a matching of size at least $k/4$ and thus Theorem \ref{fixedprob}
will imply the promised upper bound on the probability of this set 
of assignments. 

In order to do this, we need one more definition. 

\begin{definition}
Let $x,y$ be two vertices of $Z$ such that $Z$ has a path from $x$ to $y$.  
\begin{itemize}
\item  $Z(x,y)$ is the branching program obtained from $Z$ by the union of all
paths from $x$ to $y$  along with the labels on their edges.
(See Figure \ref{zxy} for an example.)
Accordingly,  for a path $P$ going through both $x$ and $y$ 
$P(x,y)$ is the subpath of $P$ starting at $x$ and ending at $y$
(again, with the labels preserved).

\item $\psi(x,y)$ is the set of all clauses $C$ of $\psi(G)$ such that for
each  path $P$ from $x$ to $y$, $A(P)$ satisfied $C$.
\item $G(x,y)$ is the subgraph of $G$ induced by the edges corresponding
to the clauses of $\psi(x,y)$.
\end{itemize}
\end{definition}

\begin{figure}[h]
\begin{tikzpicture}
\draw [fill=black]  (4,5) circle [radius=0.2];
\draw [fill=black]  (5,4) circle [radius=0.2];
\draw [fill=black]  (4,3) circle [radius=0.2];
\draw [fill=black]  (5,3) circle [radius=0.2];
\draw [fill=black]  (6,3) circle [radius=0.2];
\draw [fill=black]  (5,2) circle [radius=0.2];
\draw [fill=black]  (5,2) circle [radius=0.2];
\draw [fill=black]  (4,1) circle [radius=0.2];
\draw [fill=black]  (3,2) circle [radius=0.2];

\draw [->] (4,5) --(4.85,4.15);
\draw [->] (4,5) --(3,2.2);
\draw [->] (5,4) --(4.15,3.15);
\draw [->] (5,4) --(5,3.2);
\draw [->] (5,4) --(5.85,3.15);
\draw [->] (5,3) --(5,2.2);
\draw [->] (4,3) --(4.85,2.15);
\draw [->] (6,3) --(5.15,2.15);
\draw [->] (5,2) --(4.15,1.15);
\draw [->] (4,3) --(3.15,2.15);
\draw [->] (3,2) --(3.85,1.15);

\node [above] at (5,4.1) {$x$};
\node [below] at (5,1.9) {$y$};
\node [left] at (4.5,3.5) {$v_1$};
\node [left] at (5.1,3.5) {$v_1$};
\node [right] at (5.5,3.5) {$v_2$};
\node [left] at (4.5,2.5) {$v_2$};
\node [left] at (5.1,2.5) {$v_3$};
\node [right] at (5.5,2.5) {$v_3$};
\node [left] at (3.5,3.5) {$v_1$};
\node [left] at (3.5,1.5) {$v_4$};

\draw [dashed] (3.8,4.2) rectangle (6.1,1.8);
\end{tikzpicture}
\caption{A braching program $B$. The dashed rectangle denotes $B(x,y)$}
\label{zxy}
\end{figure}
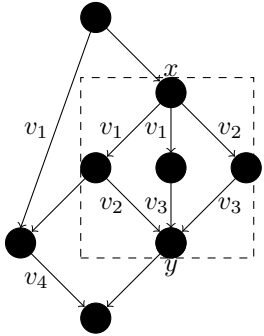
\begin{lemma} \label{broadchunk}
Let $P$ be a source-sink path of  $Z$ and let $u_1, \dots, u_{a+1}$ be yardsticks
of $P$ (recall that $a \leq d$). Then there is $1 \leq i \leq a$ such that
$G(u_{i},u_{i+1})$ is of pathwidth at least $k$.
\end{lemma}

{\bf Proof.}
Let $\psi^*=\psi(u_1,u_2) \cup \dots \cup \psi(u_{a},u_{a+1})$.
We claim that $\psi^*=\psi(G)$.
By construction, $\psi^* \subseteq \psi(G)$.
Assume that there is a $C \in \psi(G) \setminus \psi^*$.
This means that for each 
$1 \leq i \leq a$, there is a path $P_i$ from $u_{i}$ to $u_{i+1}$
such that $A(P_i)$ does not satisfy $C$. 
However, $P=P_1+ \dots+P_{a+1}$ is a source-sink path and hence
$A(P)=A(P_1) \cup \dots A(P_{a+1})$ must satisfy $\psi(G)$ while not satisfying
$C$ at the same time, a contradiction. \footnote{This is the argument where 
the monotonicity is essential. Otherwise $P_1 \cup \dots P_{a+1}$ may be an inconsistent
path.}

It follows that $G(u_1,u_2) \cup \dots \cup G(u_{a},u_{a+1})=G$.
By definition of $d$-pathwidth, the pathwidth of one of these graphs
must be at least $k$.
$\blacksquare$

In order to proceed, 
for a sequence $L$ of vertices of $Z$,
we denote by ${\bf S}(L)$ the set of satisfying assignments
of $\psi(G)$ carried by source-sink paths of $Z$ going through 
$L$. In other words, ${\bf S}(L)$ consists of all satisfying assignments
$S$ of $\psi(G)$ such that there is a source-sink path $Q$
of $Z$ going through all the vertices of $L$ such that $A(Q) \subseteq S$. 

%we denote by ${\bf P}(x,a,y)$ the set 
%of all paths of $Z(x,y)$ going through $a$.
%Accordingly, $A({\bf P}(x,a,y)}$ are the respective assignments
%carried by these paths. 

\begin{lemma}  \label{largesep}
Let $x,y$ be vertices of $Z$ such that $Z$ has a path from $x$ to $y$ and 
$Z(x,y)$ is read-once. 
Let $r=pw(G(x,y))$. Then each $x,y$-path $P$ of $Z$ has a 
node $a$ such that ${\bf S}(L)$ fixes a matching of clauses
of size at least $\lfloor r/4 \rfloor$ where $L=(x,a,y)$.
\end{lemma}

%Combination of Lemmas \ref{broadchunk} abd \ref{largesep}
%shows that any source-sink path $P$ of $Z$ has consecutive
%yardsticks $x$ and $y$ and a vertex $a$ of $P(x,y)$ such that
%${\bf P}(x,a,y)$  fixes a matching of size $k/4$.
%Since the $x,y$ subpath of each source-sink path going through $x,a,y$ 
%is an element of ${\bf P}(x,a,y)$, the assignments carried by these paths
%fix the very same matching. Thus $x,a,y$ is the desired triple of vertices
%for $P$.
Now, we are ready to provide a formal proof of Theorem \ref{mainlower}, where
the desired triple of vertices for each source-sink path $P$ of $Z$
readily follows from combination of Lemma \ref{broadchunk}
and Lemma \ref{largesep}.

{\bf Proof of Theorem \ref{mainlower}.}
In light of Theorem \ref{yardsimulate}, we may assume that
$Z$ is a $d$-SMNBP with yardsticks.

Let $P$ be a source-sink path of $Z$.
By Lemma \ref{broadchunk}, there are
consecutive yardsticks $x,y$ of $P$ such
that $G(x,y)$ is of pathwdith at least $k$.
Further on, by Lemma \ref{largesep},
there is a vertex $a$  such that 
${\bf S}(L)$ fixes a matching of clauses
of size at least $\lfloor k/4 \rfloor$, where $L=(x,a,y)$.  
By Theorem \ref{fixedprob},

\begin{equation} \label{maineq1}
Pr({\bf S}(L)) \leq 7/8^{\lfloor k/4 \rfloor}
\end{equation}

We call $L$ as above a \emph{path triple} (for $P$ if a path needs to be specified). 
Let ${\bf L}$ be the set of all path triples.
Observe that $Pr(\bigcup_{L \in {\bf L}} {\bf S}(L))=1$
Indeed, let $S$ be a satisfying assignment of $\psi(G)$.
Let $P$ be a source-sink path of $Z$ such that $A(P) \subseteq S$.
Then, by definition, $S \in {\bf S}(L)$ where $L$ is the path triple of $P$.
Combining the union bound with \eqref{maineq1}, we obtain,
that $1 \leq \sum_{L \in {\bf L}} Pr({\bf S}(L)) \leq |{\bf L}|*(7/8)^{\lfloor k/4 \rfloor}$
from where we conclude that  
\begin{equation} \label{maineq2}
|{\bf L}| \geq (8/7)^{\lfloor k/4 \rfloor}
\end{equation} 

On the other hand, by construction, 
$|{\bf L}| \leq |V(Z)|^3$ hence $|V(Z)| \geq \Omega((8/7)^{k/12})$.
As $Z$ is connected , $|Z| \geq |V(Z)|-1$ and hence the
statement follows.
$\blacksquare$

Let us discuss two directions of further research.
\begin{question}
Does Theorem \ref{mainlower} hold for $d$-MNBP (that is, without the separability assumption)?
\end{question}

Since Lemma \ref{broadchunk} is not true without the assumption that
$Z$ is monotone, the above argument does not work for the non-monotone case.
We believe that resolving the following open question will provide an important
insight in this direction.

\begin{question}
Is there a constant $\alpha$ such that for each graph $G$ of $2$-pathwidth at least $k$,
the the size of a conjunction of $2$ OBDDs representing $\psi(G)$ is at least $\alpha^k$?
\end{question}

A far fetched generalization of the last open question is whether Theorem \ref{mainlower}
holds without the monotonicity assumption but with the extra assumption that the branching
program is deterministic and oblivious. In other words, whether $d$-IBDD can be considered instead
of $d$-SMNBP.  In particular, is there an exponential lower bound for, say, $n^{0.01}$-IBDDs representing
$\psi(K_n)$?  Resolving this question positively will mean a significant breakthrough in the area of circuit 
complexity beating the $O(log^2 n)$ repetition barrier for oblivious branching programs (even though
restricted to IBDDs). On the other hand, if the question is resolved negatively, this will result in an
interesting insight as to how non-monotonicity outperforms monotonicity in this particular context.  

%\bibliographystyle{plain}
%\bibliography{KnowComp} 

\appendix
\section{Proof of Theorem \ref{yardsimulate}}
\begin{definition} \label{def:separated}
Let $Z$ be a $d$-SMNBP.
A \emph{junction vertex} of $Z$ is a vertex whose in-degree or out-degree is greater than one.
Otherwise $v$ is a \emph{subdivision vertex}. 
We say that $Z$ is \emph{subdivided} if the following two conditions hold.
\begin{enumerate}
\item {\bf Non-adjacent junctions.} There is no edge between two junction vertices.
\item {\bf Subdivided literals.} Both vertices of each edge labelled with a literal are subdivision ones.
\end{enumerate}
\end{definition}

It is not hard to see that a $d$-SMNBP can be made subdivided by subdivision of each edge into three parts
so that the label (if exists) is assigned to the middle part. This increases the size of the $d$-SMNBP
at most $3$ times. We are going to demonstrate that a subdivided $d$-SMNBP is in fact a $d$-SMNBP
with yardsticks.  In order to do this, we identify on each source-sink path a vertex we call a \emph{pre-pivot}
(the reason for this name will become clear when we provide a formal definition).
The main technical statement of the proof (Theorem \ref{pivprepiv}) demonstrates
that each path from the source to a pre-pivot is read-once and each branching program
created by the union of all paths from a pre-pivot to the sink is a $(d-1)$-SMNBP.
After that Theorem \ref{yardsimulate} is proved using a simple induction.

Throughout the proof, we often regard $Z$ as a partial order relation where $u \leq v$
if and only if $Z$ has a path from $u$ to $v$. The notions of minimal and maximal vertices
are naturally defined in this context.

\begin{definition}
A path $P$ of $Z$ is \emph{read-once} if it does not have two edges labelled by the same variable.
A vertex $v$ of $Z$ is \emph{read-once} if every path from the source to $v$ is read-once.
Otherwise $v$ is \emph{non-read-once}.  
Finally $v$ is \emph{minimally non-read-once} if every vertex $u \neq v$ such that $Z$ has a path from
$u$ to $v$ is read-once.  We denote by $M(Z)$ the set of all minimally non-read-once vertices.
\end{definition}

In the remaining part of the proof, for a vertex $v \in V(Z)$,
we denote by $Z_v$ the subgraph of $Z$ induced by $v$ and all the vertices
reachable from $v$, the labels of edges preserved. 

The following proposition is immediate by induction on the distance from the source of $Z$. 
\begin{proposition} \label{prop:mz}
For every vertex $u \in V(Z)$ that is not read-once there is $v \in M(Z)$
such that $u \in V(Z_v)$. 
\end{proposition}

\begin{definition}
Let $P$  be a source-sink path of $Z$.
Suppose that $Z$ is not read-once. 
The \emph{pivot} $w$ of $P$ is the minimal
non read-once vertex  of $P$.
The immediate predecessor of $w$ on $P$
is called the \emph{pre-pivot} of $P$.
\end{definition}

\begin{remark} \label{rem:prepivot}
\begin{enumerate}
\item Since $Z$ is not read-once, the sink of $Z$ is not read-once.
That is $P$ has non-read-once vertices. Consequently, both 
the pivot and pre-pivot of $P$ are well defined.
\item By definition, the pre-pivot of $P$ is a read-once vertex.
\end{enumerate}
\end{remark}

\begin{comment}
\begin{proposition}
Both pivot and pre-pivot of $P$ are well defined.
\end{proposition}

{\bf Proof.}
As $P$ is not read-once, the sink $sink(Z)$
is a non-read-once vertex.  Then, by definition, there is
a minimally non-read-once vertex $v$ (possibly, $sink(Z)$ 
itself) such that $sink(Z) \in V(Z_v)$. Consequently,
we can uniquely identify the minimal such a vertex $w$
on $P$, and this will be the pivot. 
Since $w$ is non-read-once $w$ is not the source and hence
has predecessors meaning that the pre-pivot is also
uniquely defined.
$\blacksquare$ 
\end{comment}

\begin{theorem} \label{pivprepiv}
Let $Z$ be a subdivided $d$-SMNBP.
Let $P$ be a non-read-once source-sink path. 
Let $v$ and $w$ be, respectively, the pre-pivot and the pivot of $P$.
Then $Z_v$  is a subdivided $d-1$-SMNBP.
\end{theorem}

{\bf Proof.}
%We start by noticing that existence of non-read-once 
%paths means that $d>1$.
The proof is divided into two cases.

{\bf Case1:  $w \in M(Z)$.}

Let $P_0$ be a non-read-once source-$w$ path.
One of the in-coming edges of $w$ 
must be labelled with a literal. 
Indeed, otherwise, we can take the predecessor of $w$ 
on $P_0$ as a non-read-once vertex in contradiction to $w \in M(Z)$.
By the second property of a subdivided $d$ -SMNBP,
$w$ is a subdivided vertex. Hence $w$ has only one in-neighbour.
As $v$ is an in-neighbour of $w$, $v$ is the only in-neigbour of $w$
and hence, due to the absence of other incoming edges, $(v,w)$ is labelled
with a literal $x$. 

%We also derive two \emph{useful properties} of $P_0$ to which we will
%refer as the \emph{first} and the \emph{second} useful property respectively.
%\begin{enumerate}
%\item The predecessor of $w$ on $P_0$ is $v$ (simply because $w$ does not have any
%other in-neighbours). 
%\item $P_0$ has a repeated occurrence of $x$ simply because otherwise $P_0 \setminus \[w\}$
%is not read-once in contradiction to the minimality of $w$.
%\end{enumerate}

Let $Q$ be a source-sink path of $Z_v$. 
By the 'Subdivided literals' property  of Definition \ref{def:separated}, 
$v$ is a subdivided vertex. Hence the first edge of $Q$ is $(v,w)$. 
Let $P^*=P_0 \setminus \{w\}+Q$. By the previous paragraph, $v$ is
the predecessor of $w$ on $P_0$ and hence
$P^*$ is a source-sink path. Hence, $P^*$ can be partitioned into 
read-once fragments $P_1, \dots P_a$ occurring in the order listed 
with $a \leq d$.  We claim that $P_1$ is a prefix of $P_0 \setminus \{w\}$
(not necessarily proper). Indeed, otherwise $P_0$ is a prefix of $P_1$
which is a contradiction as $P_0$ is not read-once.  It follows that
$Q$ is a suffix of $P_2,+ \dots+ P_a$ and hence can be partitioned into
$a-1 \leq d-1$ read-once subpaths.
Taking into account that the property of being subdivided is preserved by taking
induced subgraphs, we conclude that the theorem
holds for the considered case. 

%For the second statement, note that a non-read-once path from the source to $v$ 
%means that $v$ is a non-read-once vertex in contradiction to $w \in M(Z)$.

{\bf Case 2: $w \notin M(Z)$.}
Let $u \in M(Z)$ be such that $w \in Z_u$. 
By definition of the pivot, $v \notin Z_u$.
Consequently, in a path from $u$ to $w$ the predecessor of $w$
is not $v$, implying that the in-degree of $w$ is at least $2$ and hence
$w$ being a junction vertex. 
By the properties of Definition \ref{def:separated} , we conclude that
$(v,w)$ is not labelled and that $v$ is a subdivided vertex.

Now, let $Q$ be a source-sink path of $Z_v$.
We need to demonstrate that $Q$ can be partitioned into at most $d-1$ read-once paths.
Since $v$ is subdivided, the first edge of $Q$ is $(v,w)$.
Let $Q_w$ be the suffix of $Q$ starting at $w$.
Since $(v,w)$ is unlabelled, it is enough to show that
$Q_w$ can be partitioned into at most $d-1$ read-once paths. 
As $w$ is a pivot, there is a source-$w$ non-read-once path 
$P'_0$. It is not hard to see that $P'_0+Q_w$ is a source sink path
of $Z$. Hence $P'_0+Q_w=P_1+ \dots+P_a$ for $a \leq d$ such that
for each $1 \leq i \leq a$, $P_i$ is a read-once path.
As $P'_0$ is not read-once, we conclude that $P_1$ is a prefix of $P'_0$
and hence $Q_w$ is a suffix of $P_2+ \dots, P_a$ confirming the 
theorem. 

%desired claim.
%It remains to note that as in the previous case that 
%being a subdivided branching program is preserved by taking induced 
%subgraph. Hence, we conclude that the theorem
%holds for the considered case.

%For the second statement note that a non-read-once source-$v$ path
%means that $v$ is not read-once and hence belongs to $Z_{u'}$ of some $u' \in M(Z)$
%thus contradiction $w$ being the pivot of $P$.
$\blacksquare$

{\bf Proof of Theorem \ref{yardsimulate}.}
For each edge $e$ of $Z$ introduce two new vertices that subdivide $e$ into
a directed path $e_1,e_2,e_3$. If $e$ is labelled with a variable $x$, assign $x$ to $e_2$.
The edges $e_1$ and $e_3$ are left unlabelled.  Let $Z^*$ be the resulting branching
program. It is not hard to see that $Z^*$ is a subdivided $d$-SMNBP representing the same function as $Z$.
The theorem will immediately follow from the claim below. 

\begin{claim}
For each source-sink path $P$
of $Z^*$ there is a tuple $t(P)=(v_1, \dots, v_a)$
of vertices located on $P$ in the order listed
such that $v_1$ is the source, $v_a$ is the sink,
$a \leq d+1$, and for each $1 \leq i<a$,
each path of $Z^*$ between $v_i$ and $v_{i+1}$
is read-once. 
\end{claim}

We prove the claim by induction on $d$.
For $d=1$ simply associate each source-sink
path with the source-sink pair. Do the same if
$d>1$ but all the paths are read-once.

So, we assume that $d>1$ and that $Z^*$ has
non read-once paths. Let $P$ be a source-sink
path. By assumption, the sink of $P$ is non-read-once
hence $P$ has the pivot and pre-pivot.
Let $v$ be the pre-pivot of $P$
and let $Q$ be the suffix of $P$ starting at $v$. 
Clearly, $Q$ is a source-sink path of $Z^*_v$.
By Theorem \ref{pivprepiv},
$Z^*_v$ is $d-1$-SMNBP.

Hence, by the induction assumption, there is a tuple $(v_2, \dots v_a)$
of vertices of $Q$
such that $a \leq d+1$, $v_2=v$, $v_a$ is a the sink and for
each $2 \leq i <a$ each path of $Z^*_v$ between $v_i$ and $v_{i+1}$ 
is read-once.

Let $v_1$ be the source of $Z^*$. We claim that $(v_1, \dots, v_a)$
is the desired tuple for $P$. We only need to prove that
for each $1 \leq i <a$ each path of $Z^*$ between $v_i$ and $v_{i+1}$ is 
read-once as the rest of the statements follow by construction.
For $i=1$ this follows from the definition of pre-pivot (see Remark \ref{rem:prepivot}).
For $i>1$ this follows from the previous paragraph as each
path of $Z^*$ between two vertices of $Z^*_v$ is also a path of $Z^*_v$.  
This proves the claim and the theorem.
$\blacksquare$

\section{Proof of Theorem \ref{fixedprob}}
The proof is based on the following
idea.  Let $S_1, \dots, S_m$ be sets of variables over reals.
Suppose that we want to prove that a particular value $X$
equals $\prod_{i=1}^m \sum_{a_i \in S_i} a_i$. Then this is the same
as to prove that $X$ equals $\sum_{a_1 \in S_1, \dots, a_m \in s_m} \prod_{i=1}^m a_i$:
we simply open the brackets. In terms of probabilities, this 
idea can be expressed as the following statement.

\begin{proposition} \label{waytocount}
Let ${\bf E}_1, \dots, {\bf E}_m$ be events
and assume that each ${\bf E}_i$ is the 
disjoint union of events ${\bf E}_{i,1}, \dots, {\bf E}_{i,r_i}$. 
Let $B({\bf E}_i)=\{{\bf E}_{i,1}, \dots, {\bf E}_{i,r_i}\}$.
Assume further that for each ${\bf E'}_1 \in B({\bf E}_1),
\dots, {\bf E'}_m \in B({\bf E}_m)$,
$Pr(\bigcap_{i=1}^m {\bf E'}_i)=\prod_{i=1}^m Pr({\bf E'}_i)$. 
Then $Pr( \bigcap_{i=1}^m {\bf E}_i)=\prod_{i=1}^m Pr({\bf E}_i)$.
\end{proposition}

In order to apply Proposition \ref{waytocount},
we need to extend our terminology and to prove an auxiliary
lemma that will allow us to easily calculate probabilities
of so called \emph{guarded} assignments.

Throughout this section when we refer to an assignment $S$,
we mean that $Var(S) \subseteq Var(\psi(G))$.
Also, $S$ is the disjoint union of $S_V$ and $S_E$
where $Var(S_V) \subseteq V(G)$ and $Var(S_E) \subseteq E(G)$.  

\begin{definition} 
Let $S$ be an assignment.  %to a subset of $Var(\psi(G))$.
We denote by ${\bf Ext}(S)$ the event consisting of all
the assignments that contain $S$. 
\end{definition}

We now extend the notions of free and enforced edges as in Definition \ref{def:vertedge}
to sets of literals that do not necessarily assign all of $Var(\psi(G))$. 

\begin{definition}
Let $S$ be an assignment.
\begin{itemize}
\item Let $e \in E(G)$. 
Let $u$ and $v$ be the ends of $e$.
We say that $e$ is \emph{guarded}  (by $S$)
if $u,v \in S_V$.
The set of all guarded edges is denoted by $Guarded(S)$.
In other words,
$Guarded(S)=\{e|e=\{u,v\}\in E(G),\{u,v| \subseteq S_v\}$.
 
\item 
Let $e \in E(G)$ and let $u,v$ be the ends of $e$. 
We say that $e$ is \emph{enforced} (by $S$)
if both $u$ and $v$ occur negatively in $S_V$
Otherwise, $e$  is \emph{free}.
We denote by $Enforced(S)$ and $Free(S)$ the respective sets 
of free and enforced edges.  
%and that it is \emph{free} (w.r.t. $S_V$) otherwise. 
%We say that $S$ is \emph{guarded}  if each $e \in Var(S_E)$
%is guarded w.r.t. $S_V$ and that $S$ is \emph{valid}
%if for each $e \in Var(S_E)$ that is enforded by $S_V$,
%$e$ occurs positivley in $S$. 
%We denote by $Enforced(S)$ the set of variables 
%of $S_E$ that are enforced by $S_V$ and $Free(S)=Var(S_E) \setminus Enforced(S)$.

\end{itemize}
\end{definition}

\begin{definition}
We say that an assignment $S$ is \emph{guarded} if 
$S_E \subseteq Guarded(S)$. 
We say that $S$ is \emph{valid}  if all the variables of $S_E \cap Enforced(S)$
occur positively in $S$. 
Put it differently, an assignment is valid if it does not falsify any clauses. 
\end{definition}

\begin{lemma} \label{guardvalid}
Let $S$ be a guarded and valid assignment.
Then $Pr({\bf Ext}(S))=2^{-(|S_V|+|Free(S) \cap S_E|})$. 
\end{lemma}

{\bf Proof.}
We assume first that $Var(S_V)=V(G)$. 
%We say that $e \in E(G)$ is \emph{potentially free} w.r.t. 
%$S$ if $e \notin S_E$ and $e$ is free w.r.t. $S_V$.
%Denote by $PFree(S)$ the set of edges that are 
%potentially free w.r.t. $S$. 

Consider $S^* \in {\bf Ext}(S)$.
By definition and our assumption, $Pr(S^*)=2^{-(|S_v|+|Free(S^*)|)}$.
As the number of free variables is completely determined
by the assignment to $V(G)$, we replace $Free(S^*)=Free(S)$,
that is $Pr(S^*)=2^{-(|S_V|+|Free(S)|)}$. Note that the probability of $S^*$
is completely determined by $S$, that is, all the elements
of ${\bf Ext}(S)$ have the same probability and hence
$Pr({\bf Ext}(S))=2^{-(|S_V|+|Free(S)|)}*|{\bf Ext}(S)|$.

$S^*$ is obtained  from $S$ by assigning variables of $E(G) \setminus S_E=
(Enforced(S) \setminus S_E) \cup (Free(S) \setminus S_E)$.
The elements of $Enforced(S) \setminus S_E$ must be assigned positively.
The elements of $Free(S) \setminus S_E$ can be assigned arbitrarily.
We conclude that $|{\bf Ext}(S)|=2^{|Free(S) \setminus S_E|}$.
Substituting the quantity into the formula in the end of the previous paragraph,
we  obtain the equality as required by the lemma.

Assume now that $S_V \subset V(G)$. 
Let ${\bf S^*}$ be the set of all $2^{|V(G) \setminus S_V|}$ extensions
of $S$ assigning the rest of vertex variables. As $S$ is guarded
all elements of $S$ remain valid (otherwise, an unguarded edge variable appearing
negatively, would forbid both its ends to occur negatively).
Of course, all the elements of ${\bf S^*}$ remain guarded. 
Let $S^* \in {\bf S^*}$.
By the first part of the proof,
$Pr({\bf Ext}(S^*))=2^{-(|V(G)|+|Free(S^*) \cap S^*_E|)}$.
As $S^*\setminus S$ assigns vertex variables only,
$S^*_E=S_E$. As all the variables of $S_E$ are guarded in $S$,
them being free or not is completely determinued by $S$.
Therefore, $Free(S^*) \cap S^*_E=Free(S) \cap S_E$.
In other words
$Pr({\bf Ext}(S^*))=2^{-(|V(G)|+|Free(S) \cap S_E|)}$.
Again, we see that this quantity is the same for all $S^* \in {\bf S^*}$.
As $Pr({\bf Ext}(S))$ is the disjoint union of 
${\bf Ext}(S^*)$ for $S^* \in {\bf S^*}$, we conclude that
$Pr({\bf Ext}(S)=2^{-(|V(G)|+|Free(S) \cap S_E|)}*|{\bf S^*}|=
2^{-(|V(G)|+|Free(S) \cap S_E|)}*2^{|V(G) \setminus S_V|}=
2^{-(|S_V|+|Free(S) \cap S_E|)}$ as required. 

\begin{comment}
The assignments of the edge variables that are not in $PFree(S)$ are fixed
either because they occur in $S_E$ or because their respective ends 
both occur negatively in $S$. For the variables of $PFree(S)$,
however, all $2^{|PFree(S)|}$ assignments are possible.
It follows that $|{\bf Ext}(S)|=2^{|PFRee(S)|}$.
 
The vertex assignments of all elements of ${\bf Ext}(S)$are the same and hence the free
edges are also the same in particular $Free(S) \cup PFree(S)$.
By definition of the probablity space, for each $S^* \in {\bf Ext}(S)$,
$Pr(S)=2^{-(|S_V|+|Free(S)|+|PFree(S)|)}$. Then $Pr({\bf Ext}(S))$ is obtained by multiplication
of this probabiity by $|{\bf Ext}(S)$ as in the end of the previous paragraph, 
confirming the lemma for the considered case. 

Assume now that $S$ does not assign all the  
vertex variables.  Let ${\bf S^*}$ be the set of
all the assignments $S^*$  with $Var(S^*)=Var(S) \cup V(G)$.
In other words ${\bf S^*}$ is the set ofall the assignments
extending $S$ to all the variables of $V(G)$. 
Clearly, $|{\bf S^*}|=2^{|V(G)|-|S_V|}$.  By the proven above,
for each $S^*  \in {\bf S^*}$, $Pr({\bf Ext}(S^*))=2^{-(|V(G)+|Free(S)|)}$
As ${\bf Ext}(S)$ is the disjoint union of $\{{\bf Ext}(S^*)|S^* \in {\bf S^*}\}$
we conclude that
$Pr({\bf Ext}(S))=\sum_{S^* \in {\bf S^*}} Pr({\bf Ext}(S^*))=
2^{|V(G)|-|S_V|}*2^{-(|V(G)+|Free(S)|)}$ implying the lemma.
\end{comment}
$\blacksquare$

{\bf Proof of Theorem \ref{fixedprob}.}
In order to utilize Proposition \ref{waytocount},
we need the following claim.

\begin{claim} \label{localprod}
Let $C_1, \dots, C_q$ be a matching of clauses of $\psi(G)$.
Let $S_1, \dots, S_q$ be valid assignments 
with $Var(S_i)=C_i$ for each $1 \leq i \leq q$.
Then $Pr(\bigcap_{i=1}^q {\bf Ext}(S_i))=\prod_{i=1}^q Pr({\bf Ext}(S_i))$.
\end{claim}

{\bf Proof.}
As $Var(S_1), \dots, Var(S_q)$ are all disjoint by definition
we can consider the assignment $S=S_1 \cup \dots \cup S_q$.
It is not hard to observe that ${\bf Ext}(S)=\bigcap_{i=1}^q {\bf Ext}(S_i)$.
So, we need to prove that $Pr({\bf Ext}(S))=\prod_{i=1}^q Pr({\bf Ext}(S_i))$.

It is not hard to see that $S$ is a guarded and valid assignment.
By definition, $S$ assigns $2q$ vertex variables. Let $k=|Free(S) \cap S_E|$.
Then, by Lemma \ref{guardvalid},
$Pr({\bf Ex}(S))=2^{-(2q+k)}$.
On the other hand, each $S_i$ is also a guarded and valid assignment.
By Lemma \ref{guardvalid}, ${\bf Ext}(S_i)$ is $2^{-2}$ if the edge variable
assigned by $S_i$ is enforced and $2^{-3}$ if it is free. It is not hard to
see that an edge variable assigned by $S_i$ is free for $S$ if and only if it is free for
$S_i$. Therefore, there are precisely $k$ assignments $S_i$ for which ${\bf Ext}(S_i)=2^{-3}$
is exactly $k$. That is $\prod_{i=1}^k{\bf Ext}(S_i)=2^{-(2*(q-k))}*2^{-3k}=2^{-(2q+k)}$ as required.
$\square$

Let $C$ be a clause and $C'$ be its proper subset.
Let $CL(C,C')$ be the set of all valid assignments $S$ with $Var(C)=S$
such that at least one variable of $C'$ occurs positively in $S$.
Let ${\bf Fix}(C,C')=\bigcup_{S \in CL(C,C')} {\bf Ext}(S)$.
Note that this union is disjoint. 

Let $C_1, \dots,C_q$ be a matching of clauses
and let $C'_1, \dots, C'_q$ be their respective proper subsets. 
Let ${\bf S} \subseteq {\bf SAT}(G)$ be a set of assignments
fixing each $C_i$ with $C'_i$.
Then ${\bf S} \subseteq {\bf Fix}(C_i,C'_i)$ for each $1 \leq i \leq q$.
That is, ${\bf S} \subseteq \bigcap_{i=1}^q {\bf Fix}(C_i,C'_i)$.
Thus theorem will follow from the combination of the following 
two statements.
\begin{enumerate}
\item $Pr(\bigcap_{i=1}^q {\bf Fix}(C_i,C'_i))=\prod_{i=1}^q Pr({\bf Fix}(C_i,C'_i))$.
\item $Pr({\bf Fix}(C_i,C'_i)) \leq 7/8$ for each $1 \leq i \leq q$.
\end{enumerate}

The first statement follows from the combination 
of Proposition \ref{waytocount} and Claim \ref{localprod}.
For the second statement we refer to the proof of Claim \ref{localprod}
for a fact that for any valid assignment $S_i$ with $Var(S_i)=C_i$,
$Pr({\bf Ext}(S_i)) \geq 7/8$. At least one such $S_i$ falsifies all
the variables of $C'_i$ and hence ${\bf Fix}(C_i,C'_i)$ is disjoint
with ${\bf Ext}(S_i)$ by definition. We conclude that
$Pr({\bf Fix}(C_i,C'_i)) \leq 7/8$.
$\blacksquare$

\section{Proof of Lemma \ref{largesep}}
First of all, it is convenient to extend the notion
of a set of assignments fixing a set of clauses 
to the case where the assignments are partial and
not necessarily over the same subset of variables. 

\begin{definition} \label{extfix}
Let ${\bf S}$ be a family of sets of literals over subsets of $Var(\psi(G))$
(not necessarily over the same subset). We say that ${\bf S}$ \emph{fixes}
a set $\{C_1, \dots, C_q\}$ of clauses if for each $C_i$ there is a proper non-empty
subset $C'_i$  such that each $S \in {\bf S}$ satisfies  all of $C_1, \dots, C_q$
and for each $C_i$, $S \cap C_i \subseteq C'_i$. 
We call $C'_1, \dots, C'_q$ \emph{witnessing subsets} of $C_1, \dots, C_q$, respectively.
\end{definition}

%Recall that we denoted $(x,a,y)$ by $L$. 
For a sequence $L$ of vertices of $Z$,
we denote by ${\bf P}(L)$ the set of all paths that start at the first vertex
of $L$, end at the last one, and go through all the intermediate ones.
Accordingly, $A({\bf P}(L))=\{A(P')|P' \in {\bf P}(L)\}$.

In order to prove Lemma \ref{largesep},
we show existence of a vertex $a$ such 
$A({\bf P}(L))$ fixes a large a matching of size at  
least $r/4$ where $L=(x,a,y)$ 
(as Definition \ref{extfix} enables us to do so). 
Let us see that the statement for ${\bf S}(L)$ will follow.
Indeed, let $C_1, \dots, C_q$ be the  clauses of a matching fixed by 
$A({\bf P}(L)$ and let $C'_1, \dots, C'_q$ be their respective witnessing subsets. 
Let $S^* \in {\bf S}(L)$.
This means that there is a source-sink path $P^*$ of $Z$ such that $A(P^*) \subseteq S^*$
and $P^*$ goes through $x,a,y$. Ths means that $P^*(x,a,y) \in {\bf P}(x,a,y)$ and hence
$A(P^*(x,a,y))$ has a non-empty intersection with all of $C'_i$.
As $A(P^*(x,a,y)) \subseteq A(P^*) \subseteq S^*$, the same is true regarding $S^*$.

The advantage of considering $A({\bf P}(L))$ is that the reasoning becomes 'local',
confined to $Z(x,y)$ and $\psi(x,y)$ rather than the whole $Z$ and $\psi(G)$. 

\begin{definition}
For two vertices $a_1,a_2$ of $Z$ let us denote by $V(a_1,a_2)$ the set of variables
occurring as labels of paths between from $a_1$ to $a_2$. 
Put it differently, $V(a_1,a_2)=Var(Z(a_1,a_2))$.
\end{definition}
%Then the following lemma holds.

\begin{lemma} \label{sepvars}
Let $b$ be a vertex of $Z(x,y)$, let $C$ be a clause
of $\psi(x,y)$ and let $C'=C \cap V(x,b)$.
Assume that $\emptyset \subset C' \subset C$.
Then one of the following two statements holds.
\begin{enumerate}
\item For each path $P' \in {\bf P}(x,b,y)$,  \footnote{Note a slight abuse of notation: not using extra brackets
for ${\bf P}((x,b,y))$ for the sake of better readability.} 
$A(P')$ has a non-empty intersection with $C'$.
\item For each path $P' \in {\bf P}(x,b,y)$,
$A(P')$ has a non-empty intersection with $C \setminus C'$.
\end{enumerate} 
\end{lemma}
 
{\bf Proof.}
Assume that the statement is not true.
Then there are paths $P_1$ and $P_2$ of ${\bf P}(x,b,y)$
such that $A(P_1) \cap C'=\emptyset$ and 
$A(P_2) \cap (C \setminus C')=\emptyset$.
Note that $A(P_1(x,b))$ does not intersect with $C$.
Indeed, otherwise, by definition of $b$ it can intersect
only with $C'$ which is impossible by definition of 
$P_1$. Further, on $A(P_2(b,y))$ does not intersect with $C$.
Indeed, otherwise, $A(P_2(b,y))$ contains a literal $u$ of $C'$.
However, this literal is also contained on a path from $x$ 
to $b$. Concatenating the former to the end of the latter,
we obtain a path with a double occurrence of $u$ in contradiction
to the read-onceness of $Z(x,y)$.

It follows that $P_1(x,b)+P_2(b,y)$ is an $x,y$-path whose set of labels does not intersect
with $C$ and hence does satisfy $C$. However, this is a contradiction
with our assumption that $C$ is a clause of $\psi(x,y)$.
$\blacksquare$

We utilize Lemma \ref{sepvars} for the proof of Lemma \ref{largesep}
in the following way. We demonstrate existence of a vertex $a$ 
on $P(x,y)$ such that the condition of Lemma \ref{sepvars} w.r.t. 
$V(x,a)$ is satisfied for a matching of clauses of $\psi(x,y)$ of size
at least $pw(G(x,y))/4$. Then, by Lemma \ref{sepvars}, 
${\bf P}(x,a,y)$ fixes the matching. In order to implement this plan
we need the following result that easily follows from Theorem 5 of \cite{obddtcompsys}.

\begin{theorem} \label{pwtomw}
Let $H$ be a graph and let $\pi$ be a permutation of $V(H)$.
Then there is a prefix $\pi'$ of $\pi$ such that there is a matching
of size at least $pw(H)/2$ constsing of edges with one end in $\pi'$
and the other end in $\pi \setminus \pi'$.
\end{theorem}

In order to connect Lemma \ref{sepvars} with Theorem \ref{pwtomw},
we introduce the definition of a witnessing permutation for $P(x,y)$.

\begin{definition}
Let $x_1=x, \dots, x_m=y$ be the vertices of $P(x,y)$ occurring
on $P(x,y)$ in the order listed. 
A permutation $\pi$ of $Var(\psi(x,y))$ is a \emph{witnessing permutation}
of $P(x,y)$ if $V(x,x_1) \subseteq \dots \subseteq V(x,x_m)$ are all
(elements of some) prefixes of $\pi$.
Put it differently, a witnessing permutation of $P(x,y)$ s created as
follows. Arbitrarily order $V(x,x_2)$ (as $x=x_1$, $V(x,x_1)=\emptyset)$
and let it be the initial prefix.
Then, for each $3 \leq i \leq m$, if $V(x,x_{i-1}) \subset V(x,x_i)$.
arbitrarily order $V(x,x_i) \setminus V(x,x_{i-1})$ and append the obtained
sequence to the prefix already created. Finally, the elements 
of $Var(\psi(x,y)) \setminus V(x,y)$ (if any) are appended after the elements
of $V(x,y)$ in an arbitrary  order. 
\end{definition}

As $V(G(x,y)) \subseteq Var(\psi(x,y))$, a witnessing permutation $\pi$ of $P(x,y)$
contains a permutation $\pi_G$ of $V(G(x,y))$ as a subsequence. 
Therefore, by Theorem \ref{pwtomw}, $\pi$ has a prefix $\pi'$ 'separating' a matching of clauses 
of size at least $pw(G)/2$.  If the variables of this preifx are precisely some
$V(x,x_i)$ then we are done by Lemma \ref{sepvars}. Otherwise, we consider two cases.
In the main case, $\pi'$ gets in between two consecutive prefixes $V(x,x_i)$ and $V(x,x_{i+1})$ and
we demonstrate that one of them separates at least half of the matching
separated by $\pi'$. A formal description of this reasoning is provided in the
proof below. 
 
{\bf Proof of Lemma \ref{largesep}.}
We assume that $r$ is a multiple of $4$.
Otherwise, we adjust the value of $r$ by subtracting at most $3$
from it. %Modify the statement!!!!!

Let $x_1=x, \dots, x_m=y$ be the vertices of $P(x,y)$ occurring
on $P(x,y)$ in the order listed. 
Let $\pi$ be a witnessing permutation of $P(x,y)$.
Let $\pi_G$ be the permutation of $V(G(x,y))$ induced by $\pi$
(that is, for $u_1,u_2 \in V(G(x,y))$ $u_1$ precedes $u_2$
in $\pi_G$ if and only if $u_1$ precedes $u_2$ in $\pi$).
Let $\pi'_G$ be a prefix of $\pi_G$ such that there is a matching
$M^*_0$ of $G(x,y)$ of size at least $r/2$ such that each edge of $M^*_0$
has one end in $\pi'_G$ and one end in its complement: the existence
of such a prefix is guaranteed by Theorem \ref{pwtomw}.
Let $\pi'$ be the prefix of $\pi'_G$ having the same last element as $\pi'_G$.
Then $\pi'_G$ is a subsequence of $\pi'$ and $\pi_G \setminus \pi'_G$ is a subsequence
of $\pi \setminus \pi'$.  Let $M_0=\{(u \vee e \vee v)| \{u,v\} \in M^*_0\}$. It follows that each
clause of $M_0$ has a non-empty intersection with both $\pi'$ and $\pi \setminus \pi'$.
If the set of elements of $\pi'$ is some $V(x_i)$, we are done by Lemma \ref{sepvars}.

Otherwise,  we consider two cases.
The first, some pathological case is that $V(x,x_m) \subset \pi'$.
In other words, $\pi'$ contains all of $Var(Z(x,y))$ and some extra variables.
As $P(x,y)$ satisfies all of $\psi(x,y)$, $V(x,y)$ intersects with all the clauses
of $M_0$. As $V(x,y) \subset \pi'$, for each $C \in M_0$, $\emptyset \subset C \cap V(x,y) \subset C$.
Hence, the result follows by Lemma \ref{sepvars}.

In the second and the more intersting case, 
there is $1 \leq i <m$ such that $V(x,x_i) \subset \pi' \subset V(x,x_{i+1}$.
Let $M_1$ be the subset of $M$ consisting of $C$ such that $V(x,x_i) \cap C \neq \emptyset$.
Note that by the choice of $x_i$, $V(x,x_i) \subset \pi' \cap C \subset C$.
Therefore,  if $|M_1| \geq r/4$, we are done by Lemma \ref{sepvars}.

Otherwise, 
$|M_1| \leq r/4-1$. 
Note that the label on the edge $(x_i,x_{i+1})$  (if exists) belongs to
at most one clause $C$ of $M_0$ (because the clauses of $M_0$ are pairwise disjoint).
Let $M'_1=M_1 \cup \{C\}$ if such a $C$ exists. Otherwise, let $M'_1=M_1$.
Clearly, $|M'_1|  \leq r/4$. Let $M_2=M_0 \setminus M'_1$. Then $|M_2| \geq r/4$.
We are going to show that for each $C \in M_2$,
$\emptyset \subset V(x,x_{i+1}) \cap C \subset C$. Then the considered lemma will
immediately follow from Lemma \ref{sepvars}.

Since $\pi' \subset V(x,x_{i+1})$ and $C \cap \pi' \neq \emptyset$,
$\emptyset \subset V(x,x_{i+1}) \cap C$. For the other containement,
observe that $A(P(x,x_{i+1})) \cap C=\emptyset$. Indeed, $V(x,x_i) \cap C=\emptyset$
by definition and, since $A(P(x,x_i)) \subseteq V(x,x_i)$, $A(P(x,x_i)) \cap C=\emptyset$.
Also, the label on $(x_i,x_{i+1})$, if any, is not in $C$ by construction.
On the other hand,  as $C$ is a clause of $\psi(x,y)$, $A(P(x,y)) \cap C \neq \emptyset$.
It follows that $A(P(x_{i+1},y)) \cap C \neq \emptyset$. Let $w \in A(P(x_{i+1},y)) \cap C$.
Then $w \notin V(x,x_{i+1})$. Indeed, otherwise, there is a path $Q$ from $x$ to $x_{i+1}$
with $w \in A(Q)$. Then $Q+P(x_{i+1},y)$ is an $x,y$ path of $Z$ with a double occurrence
of $w$ in contradiction to the read-onceness of $Z(x,y)$. We conclude that
$V(x,x_{i+1}) \cap C \subset C$ as required. $\blacksquare$ 

%\bibliographystyle{plain}
%\bibliography{KnowComp} 

\end{document}